\documentclass[12pt]{article}
\usepackage{curves}
\usepackage{amsmath, amssymb, lamsarrow}
\usepackage{amssymb}
\usepackage{amsthm}
\usepackage{graphicx}
\usepackage[english]{babel}
\def\mineappendix{
        \setcounter{section}{1}
        \setcounter{subsection}{0}
        \def\thesection{\Alph{section}}
        \def\sectionap{\@startsection  {section}{1}{\z@}
                        {-3.5ex plus-1ex minus-.2ex} {0ex plus.2ex}
                        {\reset@font\Large\bf  Appendix:  \, }
                        }
        }
\makeatother
\def\Proclaim #1. #2\par{\bigbreak\noindent{\sc#1.\enspace}{\it#2}\par}

\newcommand{\vd}{\textbf{\emph{d}}}
\newcommand{\vk}{\textbf{\emph{k}}}

\newtheorem{thm}{Theorem}[section]
\newtheorem{cor}[thm]{Corollary}
\newtheorem{lem}[thm]{Lemma}
\newtheorem{rem}[thm]{Remark}
\newtheorem{ex}[thm]{Example}
\newtheorem{prop}[thm]{Proposition}
\newtheorem{conj}[thm]{Conjecture}
\newtheorem{conv}[thm]{Convention}

\title{On A Tautological Relation Conjectured By Buryak-Shadrin}
\author{Xiaobo Liu \thanks{Research was partially supported by NSFC grants 11890662 and 12341105.}, \,\,\, \,\,\, Chongyu Wang}
\date{}

\begin{document}
\maketitle
\allowdisplaybreaks

\begin{abstract}
Buryak and Shadrin conjectured a tautological relation on moduli spaces of curves $\overline{\mathcal{M}}_{g,n}$ which has the form  $B^m_{g, \textbf{\emph{d}}}=0$  for certain tautological classes $B^m_{g, \textbf{\emph{d}}}$  where $m \geq 2, n \geq 1$ and $|\textbf{\emph{d}}| \geq 2g+m-1$. In this paper we prove that this conjecture holds if it is true for the $m=2$ and $|\textbf{\emph{d}}| = 2g+1$ case. This result reduces the proof of this conjecture to checking finitely many cases for each genus $g$. We will also prove this conjecture for the $g=1$ case.
\end{abstract}

\section{Introduction}\label{section introduction}

Let $\overline{\mathcal{M}}_{g,n}$ be the moduli space of genus $g$ stable curves with $n$ marked points.
Relations among tautological classes on $\overline{\mathcal{M}}_{g,n}$ are fascinating topics and
have been extensively studied in algebraic geometry. See, for example, \cite{P} for a nice survey on such relations.
Relations involving only $\psi$-classes and boundary classes are particularly important since they
give universal equations for generating functions of Gromov-Witten invariants (or more generally for
potential functions of cohomological field theories). A method for translating such tautological relations to
universal equations can be found in \cite{L1} and \cite{L2}. In recent years, it was also found that
such relations play important roles in the study of integrable systems associated to cohomological field theories
and their generalizations (see, for example, \cite{BHIS}, \cite{BS}, \cite{DLYZ} and \cite{HIS}).
The purpose of this paper is to study a series of tautological relations conjectured by Buryak and Shadrin in \cite{BS}.

In \cite{BS}, the authors introduced a set of tautological classes
$B^m_{g, \textbf{\emph{d}}}$ which can be roughly defined as follows.
Let $T$ be a genus-$g$ stable rooted tree with $n$ regular legs, $m$ frozen legs, and $k$ extra legs.
$T$ can be considered as a dual graph of certain stable curves and it defines a boundary strata
 of ${\overline {\cal M} }_{g, n+m+k}$.
Let $\textbf{\emph{d}}=(d_1, \cdots, d_n)$ be an $n$-tuple of nonnegative integers. Together
with the distribution of extra legs, $\textbf{\emph{d}}$ induces a $\psi$-decoration $q_{\textbf{\emph{d}}}$
on this boundary strata which gives the powers of $\psi$-classes associated to half edges of $T$.
One can define a tautological class
$e_* [T , q_{\textbf{\emph{d}}}] \in H^{2|\textbf{\emph{d}}|} ({ \overline {\cal M} }_{g, n+m})$
where $e$ is the map forgetting all marked points labelled by extra legs  and $|\textbf{\emph{d}}|:= d_1 + \cdots + d_n$.
Then
\[B^m_{g, \textbf{\emph{d}}} :=  \sum_{T \in {\textrm{SRT}}^{(b,nd)}_{g,n,m;o}} (-1)^{|E(T)|} e_* [T , q_{\textbf{\emph{d}}}] \in H^{2|\textbf{\emph{d}}|} ({ \overline {\cal M} }_{g, n+m}),\]
where
${\textrm{SRT}}^{(b,nd)}_{g,n,m;o}$ is the set of nondegenerate balanced stable genus $g$ rooted trees with $m$ frozen legs and $n$ regular legs, $|E(T)|$ is the number of edges on $T$.
The precise definition of $B^m_{g, \textbf{\emph{d}}}$ and terminologies involved here will be given
in Section \ref{sec:BSc}.

Buryak and Shadrin conjectured that
$B^m_{g, \textbf{d}}=0$
if $m \geq 2$ and $|\textbf{d}| \geq 2g+m-1$.
They have showed that this conjecture can be applied to prove the polynomiality of the Dubrovin-Zhang hierarchy
associated with an arbitrary F-Cohomological Field Theory. On the other hand,
since $B^m_{g, \textbf{d}}$ only involves $\psi$-classes and boundary classes, this conjecture also gives
a series of universal equations for Gromov-Witten invariants.
It was also proved in \cite{BS} that this conjecture holds if either $n=1$ or $g=0$.
Moreover, it was showed that the proof of this conjecture can be reduced to checking the case
$|\textbf{d}| = 2g+m-1$ with all components of $\vd$ positive.
The main result of this paper is that we can further reduce the proof of this conjecture to the case where $m=2$, i.e.
we have the following
\begin{thm} \label{thm:m=2}
If $B^2_{g, \textbf{d}}=0$ for all $g \geq 0$, $n \geq 1$, and $\vd \in \mathbb{Z}^n$ with all components of $\vd$ positive
and $|\textbf{d}| = 2g+1$,
then  $B^m_{g, \textbf{d}}=0$ for all $m \geq 2$, $g \geq 0$, $n \geq 1$,
and $\vd \in (\mathbb{Z}_{\geq 0})^n$ with $|\textbf{d}| \geq 2g+m-1$.
\end{thm}
\noindent
This result reduces the proof of Buryak-Shadrin conjecture to checking a finite number (i.e. the number of partitions of $2g+1$) of equations for each fixed genus $g$.
In particular, we will prove that this conjecture holds for the $g=1$ case.
\begin{thm}\label{main theorem}
For any $m \geq 2,  n \geq 1$, and $n$-tuple of nonnegative integers $\textbf{d}$. If $|\textbf{d}| \geq m+1$, then
$B^m_{1, \textbf{d}} =0$.
\end{thm}
We believe the same method also works for the genus 2 case, which will be studied in a separate paper.

This paper is organised as follows. In Section \ref{section preliminary}, we review definition and
basic properties of tautological classes, set up basic notations, and give precise definition
for the class $B^m_{g, \textbf{d}}$.
 In Section \ref{section push of B}, we prove a push-forward formula for $B^m_{g, \textbf{\emph{d}}}$, which will be used in Section \ref{section induction of thm} to give a proof of Theorem \ref{thm:m=2}.
 Proof for Theorem \ref{main theorem}
 will be given in Section~\ref{section computation of base case}.

\section{Prelimilary}\label{section preliminary}

In this section, we will recall some basic definitions for moduli spaces of curves and dual graphs
and set up basic notations used in this paper. We will follow the terminologies used in
\cite{ACG} and \cite{BS}.

\subsection{Dual graphs and tautological classes}

Recall that a {\it dual graph} $\Gamma$ consists of the following data:

\begin{itemize}

\item A finite nonempty set of vertices $V(\Gamma)$;

\item A finite set of half edges $H(\Gamma)$;

\item An involution map $\iota: H(\Gamma) \longrightarrow H(\Gamma)$ with $\iota^2=id$;

\item An attachment for each half edge $h \in H(\Gamma)$ a unique vertex $v(h) \in V(\Gamma)$;

\item An assignment for each vertex $v$ a non-negative integer $g_v$, called its genus.

\end{itemize}

\noindent
Fixed points of $\iota$ are called {\it legs} of $\Gamma$. The set of all legs of $\Gamma$ is denoted by $L(\Gamma)$.
Each orbit of $\iota$ consisting of two distinct half edges is called an {\it edge} of $\Gamma$.
The set of all edges of $\Gamma$ is denoted by $E(\Gamma)$.
For each vertex $v$, the set of all half edges attached to $v$ is denoted by $H_v$.
A dual graph $\Gamma$ is {\it stable} if $2g_v -2+|H_v|>0$ for all $v \in \Gamma$, where $|H_v|$ is
the number of elements in $H_v$.
The {\it genus of a dual graph} $\Gamma$ is defined to be
\[ g(\Gamma) := 1+|E(\Gamma)|-|V(\Gamma)| + \sum_{v \in V(\Gamma)} g_v.\]
Unless otherwise stated, we will always assume $\Gamma$ is connected.

Let $\overline{\cal M}_{g, n}$ be the moduli space of genus $g$ stable nodal curves with $n$ marked points.
Given any $C \in \overline{\cal M}_{g, n}$, one can associate a dual graph whose vertices correspond to irreducible components of $C$, edges correspond to nodal points on $C$, and legs correspond to marked points on $C$.
On the other hand, given a stable dual graph $\Gamma$, there is a {\it clutching map}
\[ \xi_\Gamma : \prod_{v\in V(\Gamma)} \overline{\mathcal{M}}_{g_v , |H_v|} \to \overline{\mathcal{M}}_{g(\Gamma), |L^(\Gamma)|} \]
defined in the following way.
Given $C_v \in \overline{\mathcal{M}}_{g_v , |H_v|}$ for $v\in V(\Gamma)$, the curve $\xi_\Gamma( \prod_{v\in V(\Gamma)} C_v)$ is obtained by gluing curves $C_v$ along certain marked points according to the following rule:
If $h \in H_v$ is not a leg on $\Gamma$, then the marked point on $C_v$ corresponding to $h \in H_v$
should be glued to the marked point on $C_{v^*}$ corresponding to $\iota(h)$,  where $v^*=v(\iota(h))$.
If $C_v$ is smooth for all $v \in V(\Gamma)$, then $\xi_\Gamma( \prod_{v\in V(\Gamma)} C_v)$ is a curve with dual graph isomorphic to $\Gamma$.

Fix a label $i$ for marked points. There is a line bundle ${\mathbb L}_i$ over $\overline{\cal M}_{g, n}$
whose fibre over a curve $C \in \overline{\cal M}_{g, n}$ is the cotangent space of $C$ at the marked point labelled by $i$.
The first Chern class of ${\mathbb L}_i$ is called the $\psi$-class associated to the marked points labelled by $i$ and is denoted by $\psi_i$.

Given a stable dual graph $\Gamma$ and $v \in V(\Gamma)$, we can use $H_v$ as the set of labels for marked points on
curves $C_v \in \overline{\mathcal{M}}_{g_v , |H_v|}$. So for each $h \in H_v$, we have a $\psi$-class, denoted by $\psi_h$, on $\overline{\mathcal{M}}_{g_v , |H_v|}$
associated to the marked point labelled by $h$. The following simple observation will be useful in calculations:
If $h$ is a leg in a dual graph $\Gamma$, it can be viewed as a label for marked points on curves in both spaces $\overline{\mathcal{M}}_{g_{v(h)} , |H_{v(h)}|}$ and $\overline{\mathcal{M}}_{g(\Gamma), |L^(\Gamma)|}$. Therefore it defines
 a class $\psi_h$ on both spaces. Let
\[ P_h:  \prod_{v\in V(\Gamma)} \overline{\mathcal{M}}_{g_v , |H_v|} \to \overline{\mathcal{M}}_{g_{v(h)} , |H_{v(h)}|} \]
be the natural projection map. Then it is well known that
\[  {P_h}^* (\psi_h)  = {\xi_\Gamma}^* (\psi_h), \]
where $\psi_h$ on the left hand side is a cohomology class on  $\overline{\mathcal{M}}_{g_{v(h)} , |H_{v(h)}|}$
and $\psi_h$ on the right hand side is a cohomology class on $\overline{\mathcal{M}}_{g(\Gamma), |L^(\Gamma)|}$
(see, for example, equation (4.30) on page 581 in \cite{ACG}).
Hence by the push-pull formula, we have
\begin{equation}\label{11111}
 {\xi_\Gamma}_* ({P_h}^* (\psi_h) \cdot \alpha) = \psi_h \cdot {\xi_\Gamma}_* (\alpha),
\end{equation}
for any cohomology class
$\alpha \in H^*(\prod_{v\in V(\Gamma)} \overline{\mathcal{M}}_{g_v , |H_v|})$.
Note that equation \eqref{11111} may not apply if $h$ is not a leg.

Given a stable dual graph $\Gamma$, a {\it $\psi$-decoration} on $\Gamma$ is defined to be a map
$$q : H(\Gamma) \to \mathbb{Z}_{\geq 0}.$$
The pair $(\Gamma , q)$ is called  a {\it $\psi$-decorated graph}. It gives a
tautological class
\begin{equation}\label{!!!}
[\Gamma , q] : = \frac{1}{|Aut(\Gamma)|} \,\, {\xi_\Gamma}_* \bigg(\bigotimes_{v\in V(\Gamma)} \,\, \prod_{h\in H_v} \psi^{q(h)}_h \bigg)
        \in H^* (\overline{\mathcal{M}}_{g(\Gamma), |L(\Gamma)|}),
\end{equation}
where $Aut(\Gamma)$ is the group of automorphisms of $\Gamma$.
For the convenience of calculations and typesetting, we introduce the following notation
to represent the tautological class $[\Gamma , q]$:
\begin{equation} \label{eqn:taut}
\prod_{v\in V(\Gamma)}<\prod_{h\in H_v} \Psi^{q(h)}(h)>_{g_v}.
\end{equation}
In such expressions, we will also use $h^*$ to represent $\iota(h)$  if $h \in H(\Gamma)$ is not a leg. Moreover we will identify
$\Psi^{0}(h)$ and $h$. Note that expression \eqref{eqn:taut} contains all information of the
dual graph $\Gamma$ and the $\psi$-decoration $q$.
In fact, given an expression of the form \eqref{eqn:taut} where $V(\Gamma)$ is replaced by any set $V$
and $H_v$ just means a set indexed by $v$, we can find a unique dual graph $\Gamma$ and
a $\psi$-decoration $q$ on $\Gamma$ such that the tautological class $[\Gamma, q]$ is represented
by expression \eqref{eqn:taut}. We will illustrate this procedure by the following example.

\begin{ex}
Let $(\Gamma, q)$ be a $\psi$-decorated graph given by the following picture
\begin{center}
\includegraphics[scale=0.55]{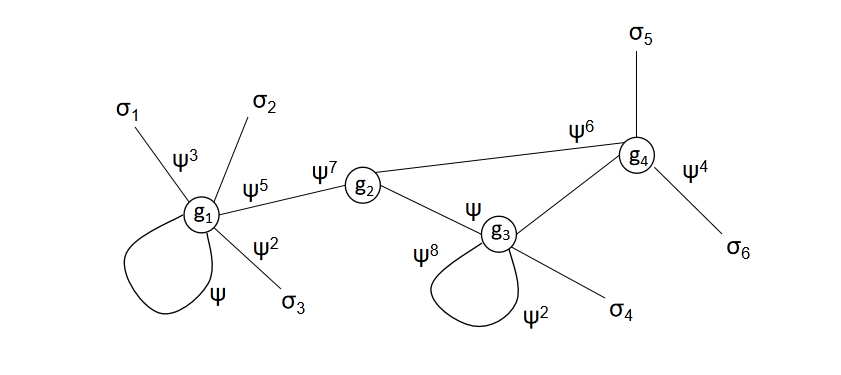}
\end{center}
where $\sigma_1, \cdots, \sigma_6$ represent legs of $\Gamma$,
each vertex $v_i$ with genus $g_i$ is represented by a circle around the number $g_i$, and
the $\psi$-decoration is represented by the powers of $\psi$ nearby the corresponding half edges.
Assume edges of $\Gamma$ are labelled in the following way
\begin{itemize}
\item $\{ \gamma_1, \gamma_1^*\}$ is the edge connecting vertex $v_1$ to itself,
\item $\{ \gamma_2, \gamma_2^*\}$  is the edge connecting vertices $v_1$ and $v_2$,
\item $\{ \gamma_3, \gamma_3^*\}$  is the edge connecting vertices $v_2$ and $v_3$,
\item $\{ \gamma_4, \gamma_4^*\}$  is the edge connecting vertices $v_2$ and $v_4$,
\item $\{ \gamma_5, \gamma_5^*\}$  is the edge connecting vertex $v_3$ to itself,
\item $\{ \gamma_6, \gamma_6^*\}$  is the edge connecting vertices $v_3$ and $v_4$.
\end{itemize}
Then the tautological class given by the above $\psi$-decorated graph can be represented as
\begin{eqnarray}
&&<\Psi^3 (\sigma_1) \sigma_2 \Psi^2 (\sigma_3) \gamma_1 \Psi(\gamma^*_1) \Psi^5 (\gamma_2)>_{g_1}
<\Psi^7 (\gamma^*_2) \gamma_3 \gamma_4>_{g_2} \nonumber \\
&& \hspace{60pt} <\Psi (\gamma^*_3) \Psi^8 (\gamma_5) \Psi^2 (\gamma_5^*) \sigma_4 \gamma_6>_{g_3}
<\Psi^6 (\gamma^*_4) \gamma^*_6 \sigma_5 \Psi^4 (\sigma_6)>_{g_4}.  \label{eqn:tautex}
\end{eqnarray}

On the other hand, given expression \eqref{eqn:tautex}, we can easily recover the dual graph $\Gamma$ in the following way:
Each angle bracket $< \cdots >_{g_i}$ gives a vertex $v_i$ with genus $g_i$. Each $\Psi^q(h)$ (or simply $h$ if $q=0$)
inside the angle bracket $< \cdots >_{g_i}$ represents a half edge (still denoted by $h$) attached to the vertex $v_i$. Any pair of
distinct half edges of the form $\{h, h^*\}$ gives an edge connecting the vertices attached to $h$ and $h^*$ respectively.
A leg is an half edge $h$ for
which $h^*$ does not appear in expression \eqref{eqn:tautex}. We can also recover the $\psi$-decoration
$q: H(\Gamma) \longrightarrow \mathbb{Z}_{\geq 0}$ by letting $q(h)=q$ if  $\Psi^q(h)$ appears in expression \eqref{eqn:tautex}. This procedure can be easily generalized for any expression of the form \eqref{eqn:taut}.

We should point out that relabelling edges does not change the dual graph. So we could freely change
labels of edges if it is needed. For example, we can also choose $\{ \gamma_2, \gamma_2^* \}$ to be the
edge connecting vertex $v_1$ to itself, and choose  $\{ \gamma_1, \gamma_1^* \}$ to be the edge connecting
vertices $v_1$ and $v_2$. Then the tautological class represented by expression \eqref{eqn:tautex} can also
be represented by
\begin{eqnarray}
&&<\Psi^3 (\sigma_1) \sigma_2 \Psi^2 (\sigma_3) \gamma_2 \Psi(\gamma^*_2) \Psi^5 (\gamma_1)>_{g_1}
<\Psi^7 (\gamma^*_1) \gamma_3 \gamma_4>_{g_2} \nonumber \\
&& \hspace{60pt} <\Psi (\gamma^*_3) \Psi^8 (\gamma_5) \Psi^2 (\gamma_5^*) \sigma_4 \gamma_6>_{g_3}
<\Psi^6 (\gamma^*_4) \gamma^*_6 \sigma_5 \Psi^4 (\sigma_6)>_{g_4}.
\end{eqnarray}

\end{ex}

It is also convenient to allow $q(h) < 0$ for some $h$ in the expression \eqref{eqn:taut}. In this case, the corresponding tautological class is just understood as $0$.

\subsection{Buryak-Shadrin Conjecture}
\label{sec:BSc}

 Recall that a tree $T$ is a connected dual graph whose first Betti number is $0$.
 In a tree $T$, any two vertices $v_1, v_2 \in V(T)$ are connected by a unique path
 such that an edge can occur in this path at most once. This path is denoted by
 $\overline{v_1 v_2}$.

 A {\it rooted tree} is a tree
 with a distinguished vertex, called the {\it root}. We usually denote the root by $v_0$ unless otherwise stated.
 On a rooted tree $T$, we can define the {\it level} of a vertex $v \in V(T)$, denoted by $l(v)$,
 to be the number of edges in the path $\overline{v_0 \, v}$ plus $1$. In particular $l(v_0)=1$.
 Note that this is the canonical level function used in \cite{BS}.
 We call $v \in V(T)$ a {\it top vertex} if it does not lie in any path $\overline{v_0 \, v'}$ with
 $v \neq v' \in V(T)$.

As in \cite{BS}, we will consider rooted trees $T$ whose set of legs is a disjoint union of three subsets, i.e.
 \[ L(T) = L^- (T)  \coprod L^+ (T) \coprod L^0 (T), \]
 where elements of $L^- (T)$ are called {\it frozen legs}, elements of $L^+(T)$ are called {\it regular legs}, and
 elements of $L^0(T)$ are called {\it extra legs}. We will assume all frozen legs are attached to the root.
 A leg is {\it positively directed} if it is a regular leg, {\it negatively directed} if it is a frozen leg,
 {\it not directed} if it is an extra leg.
 If a half edge $h \in H(T)$ is not a leg, then $h$ and $h^*=\iota (h)$ are attached to two different vertices $v(h)$ and $v(h^*)$
 respectively. We say that $h$ is {\it positively directed} if the level of $v(h)$ is strictly smaller than the level of $v(h^*)$.
 Otherwise we say that $h$ is {\it negatively directed}.
Let $H^+ (T)$ be the set of all positively directed half edges and $H^- (T)$ the set of all negatively directed half edges. Then
 \[ H(T) = H^- (T) \coprod L^0(T) \coprod H^+ (T).\]

 A rooted tree is {\it balanced} if there are no extra legs on the root and any other vertex must have
at least one extra leg attached to it.
Given a balanced stable rooted tree $T$ with $n \geq 1$ regular legs $\{r_1, \cdots, r_n\}$ and an $n$-tuple of non-negative integers
\[ \textbf{\emph{d}} = (d_1, d_2 , \ldots , d_n) \in ({\mathbb Z}_{\geq 0})^n, \]
 we can define a $\psi$-decoration
$ q_{\textbf{\emph{d}}} : H(T)\to \mathbb{Z}_{\geq 0} $
such that
\begin{equation} \label{eqn:qd}
q_{\textbf{\emph{d}}}(h) := \left\{ \begin{array}{ll} 0 &  {\rm if} \,\,\, h \in H^-(T) \bigcup L^0(T), \\
                d_i & {\rm if} \,\,\, h=r_i \in L^+(T), \\
            |H_{v(\iota (h))}\cap L^0 (T)|-1 & {\rm if} \,\,\, h\in H^+(T)\setminus L^+(T).
        \end{array} \right.
\end{equation}
Note that $|H_v \cap L^0 (T)|$ is the number of extra legs on the vertex $v$.
 We call $d_i$ the {\it weight} of $r_i \in L^+(T)$ and $\textbf{\emph{d}}$ the weight of $T$. We then have a tautological class
 \[ [T , q_{\textbf{\emph{d}}}] = {\xi_T}_* \bigg(\bigotimes_{v\in V(T)} \prod_{h\in H_v} \psi^{q_{\textbf{\emph{d}}}(h)}_h \bigg)
        \in H^* (\overline{\mathcal{M}}_{g(T), |L(T)|}). \]
 Note that $|Aut(T)| = 1$ in this case.

A rooted tree is {\it nondegenerate} if it is still stable after forgetting all extra legs.
Let ${\textrm{SRT}}^{(b,nd)}_{g,n,m;o}$ be the set of all nondegenerate balanced genus-$g$ stable rooted trees with
$n$ regular legs, $m$ frozen legs, and possibly some extra legs.
Given $T \in {\textrm{SRT}}^{(b,nd)}_{g,n,m;o}$, we can use the set $L(T)$ to label marked points on curves in $\overline{\cal M}_{g, |L(T)|} $
and define the map
\begin{equation} \label{eqn:e}
 e : \overline{\cal M}_{g, |L(T)|} \to \overline{\cal M}_{g, n+m},
\end{equation}
which forgets all marked points corresponding to extra legs.
When we want to emphasize the dependence of $e$ on $T$, we also write $e$ as $e_T$.
Given integers $g, m \geq 0$, $n \geq 1$, and
$\textbf{\emph{d}}=(d_1, \cdots, d_n) \in ({\mathbb Z}_{\geq 0})^n$,
 define the tautological class
\begin{equation} \label{eqn:Bmgdold}
 B^m_{g, \textbf{\emph{d}}} :=  \sum_{T \in {\textrm{SRT}}^{(b,nd)}_{g,n,m;o}} (-1)^{|E(T)|} e_* [T , q_{\textbf{\emph{d}}}]
   \,\,\, \in H^{2|\textbf{\emph{d}}|} ({ \overline {\cal M} }_{g, n+m}),
 \end{equation}
where $|\textbf{\emph{d}}|:=d_1+ \cdots + d_n$.
In this equation, we consider trees obtained from $T$ by permuting labels of extra legs to be the same tree as $T$,
and therefore only appear once in the summation on the right hand side of this equation.
By Proposition 3.6 in \cite{BS}, Conjecture 1 in \cite{BS} is equivalent to the following
\begin{conj}[Buryak-Shadrin]
\label{conj}
For any $g \geq 0, m \geq 2, n\geq 1$ and $n$-tuple of nonnegative integers $\textbf{d}$ such that $|\textbf{d}| \geq 2g+m-1$, we have $B^m_{g, \textbf{d}}=0$.
\end{conj}

\noindent
In \cite{BS}, this conjecture was proved for the case where either $g=0$ or $n=1$ (using a relation
proved in \cite{LP}). Theorem \ref{main theorem} in this paper
just says that
this conjecture also holds for $g=1$.

\begin{rem}
The class $B^m_{g, \textbf{d}}$ defined by equation \eqref{eqn:Bmgdold} is in fact the class denoted by
$\Tilde{B}^m_{g, \textbf{\emph{d}}}$ in \cite{BS}.
In \cite{BS}, the notation $B^m_{g, \textbf{d}}$ was used to represent a slightly different class
which will not appear in this
paper. We borrow the same notation here for simplicity.
\end{rem}

To compute $e_* [T , q_{\textbf{\emph{d}}}]$ in the definition of $B^m_{g, \textbf{\emph{d}}}$,
we need to use the string equation.
Given any vector $\textbf{\emph{q}} = (q_1, \cdots, q_n) \in (\mathbb{Z}_{\geq 0})^n $ and positive integer $l \leq |\textbf{\emph{q}}|$
where $|\textbf{\emph{q}}| := \sum^n_{i=1} q_i$ , define
$D_l(\textbf{\emph{q}})$ to be the set of all vectors
$\textbf{\emph{p}}=(p_1, \cdots, p_n) \in \mathbb{Z}^n$ such that
$0 \leq p_i \leq q_i$ for all $1 \leq i \leq n$ and $|\textbf{\emph{p}}|= |\textbf{\emph{q}}|-l$.
Let $\pi_{n,l} : \overline{\mathcal{M}}_{g, n+l} \to \overline{\mathcal{M}}_{g, n}$ be the map
which forgets the last $l$ market points on stable curves. Then the {\it string equation} has the following form
(see, for example, \cite{W}, \cite{ACG} and \cite{BS}):
\begin{equation}\label{string equation}
(\pi_{n,l})_* \bigg(\prod^n_{i=1} \psi^{q_i}_i \bigg) = \sum_{\textbf{\emph{p}} \in D_l(\textbf{\emph{q}})}
                  \frac{l!}{\prod^n_{i=1} (q_i -p_i)!} \prod^n_{i=1} \psi^{p_i}_i.
\end{equation}

When computing $e_* [T , q_{\textbf{\emph{d}}}]$, it is convenient to use expression \eqref{eqn:taut}
to represent the class $[T , q_{\textbf{\emph{d}}}]$.
In fact, for any $\psi$-decorated stable rooted tree $(T, q)$, we can write
\[ [T, q]=\prod_{v \in V(T)} <\prod_{h \in H_v} \Psi^{q(h)}(h)>_{g_v}. \]
Then
\begin{equation} \label{eqn:e*T}
e_* [T, q]=\prod_{v \in V(T)} e_* <\prod_{h \in H_v} \Psi^{q(h)}(h)>_{g_v},
\end{equation}
where $e_* <\prod_{h \in H_v} \Psi^{q(h)}(h)>_{g_v}$ can be computed by equation \eqref{string equation}
and is a summation of many terms after getting rid of $e_*$.
To interpret the right hand side of equation \eqref{eqn:e*T} correctly, we require that the product
$\prod_{v \in V(T)}$ satisfies a formal distribution law. After expanding this product,
the right hand side of equation \eqref{eqn:e*T} becomes a summation of some expressions of the form \eqref{eqn:taut}, and therefore
can be interpreted as a summation of corresponding tautological classes.
Note that equation \eqref{eqn:e*T} corresponds to a recursion formula in the proof of Lemma 3.8 in \cite{BS}.

Observe that there are many redundant terms in the summation on the right hand side of equation \eqref{eqn:Bmgdold}.
For example, if $T$ has a top vertex $v$ which does not have any regular legs attached to it,
then there is no positively directed half edge on $v$. So $q_{\vd}(h)=0$ for all $h \in H_v$ by definition of $q_{\vd}$.
Note that if the root of $T$ is a top vertex, then $T$ can have only one vertex and all regular legs, which
 always exist since $n \geq 1$, must lie on this vertex. Hence $v$ can not be the root and
 there exists at least one extra leg on $v$ since $T$ is balanced. So by the string equation,
 $e_* [T , q_{\textbf{\emph{d}}}]=0$.

 Motivated by the above example, we introduce the following definition:
An {\it acceptable tree} is a balanced nondegenerate stable rooted tree $T$ such that every top vertex of $T$
 is attached by at least one regular leg.
The set of all acceptable trees is denoted by $\mathcal{A}$, and we use ${\mathcal{A}}_{g,n,m}$ to denote  the set of
all genus g acceptable trees with $n$ regular legs, $m$ frozen legs, and possibly some extra legs. Then we have
\begin{equation} \label{eqn:Bmgd}
B^m_{g, \textbf{\emph{d}}} =  \sum_{T \in {\mathcal{A}}_{g,n,m}}  (-1)^{|E(T)|} e_* [T , q_{\textbf{\emph{d}}}].
\end{equation}
In the rest part of this paper, we will use this equation as the definition of  $B^m_{g, \textbf{\emph{d}}}$.

\begin{conv}\label{convention}
Unless otherwise stated,  regular legs of a tree will be represented by $U$ or $U_i$,
frozen legs will be represented by $V$ or $V_i$,
and extra legs will be represented by $W, W_i$ or $W^k_i$.  The root of a tree will be represented by $v_0$.
\end{conv}

\section{Computing push-forwards of $B^m_{g, \textbf{\emph{d}}}$ }\label{section push of B}

In this section, we will prove the following  proposition which will be needed in proving
 Theorem \ref{thm:m=2} by induction on $m$.
\begin{prop}\label{prop: cor 2}
For any integers $m \geq 2$, $n \geq 1$, and $\vd \in (\mathbb{Z}_{\geq 0})^n$, we have the following equation:
\begin{equation}\label{push forward of B 2}
(\pi^{m,n}_l)_* (B^{m+l}_{g, \textbf{d}})=  \sum_{\textbf{k} \, \in D_l(\vd)} \frac{l!}{\prod^n_{i=1} (d_i - k_i)!} B^m_{g, \textbf{k}} \, ,
\end{equation}
where  $\pi^{m,n}_l : { \overline {\cal M} }_{g, n+m+l} \to { \overline {\cal M} }_{g, n+m}$ is the map forgetting
marked points labeled by the last $l$ frozen legs, and $D_l(\vd)$ has been defined before equation \eqref{string equation}.
\end{prop}

We need some preparation before proving this proposition.
Let ${\mathcal{A'}}_{g,n,m}$ be the set of all genus-$g$ stable rooted trees which have $n$ regular legs, $m$ frozen legs, and no
extra legs such that every top vertex has at least one regular leg attached to it. There is a canonical map
\begin{equation}
e':  {\mathcal{A}}_{g,n,m} \longrightarrow {\mathcal{A'}}_{g,n,m}
\end{equation}
such that $e'(T)$ is the tree obtained from $T \in {\mathcal{A}}_{g,n,m}$ by forgetting all extra legs
on $T$.

Given any $T' \in {\mathcal{A'}}_{g,n,m}$, define
\begin{equation} \label{eqn:BT'}
 B(T', \textbf{\emph{d}}) :=  \sum_{
        \scriptstyle T \in {\mathcal{A}}(T')}
                 e_* [T , q_{\textbf{\emph{d}}}] \,\, \in H^{2 |\textbf{\emph{d}}|} ({ \overline {\cal M} }_{g, n+m}),
\end{equation}
where
\begin{equation*}
\mathcal{A}(T') := \{ T \in {\mathcal{A}}_{g,n,m} \mid e'(T)=T' \}.
\end{equation*}
Then clearly
\begin{equation} \label{eqn:B=BT'}
B^m_{g, \textbf{\emph{d}}} = \sum_{T' \in {\mathcal{A'}}_{g,n,m}} (-1)^{|E(T')|} B(T', \textbf{\emph{d}}).
\end{equation}
Note that $V(T)=V(T')$ for any $T \in \mathcal{A}(T')$. Each $T \in \mathcal{A}(T')$ determines a map
\[ p: V(T')\setminus \{v_0\} \longrightarrow \mathbb{Z}_{\geq 0}\] such that $p(v)+1$ is the number of extra legs on the vertex $v$ which
is not the root $v_0$. On the other hand,
given any map $p: V(T') \setminus \{v_0\} \longrightarrow \mathbb{Z}_{\geq 0}$, there exists a unique
tree $T_p \in \mathcal{A}(T')$ which is obtained from $T'$ by adding $p(v)+1$ extra legs on every vertex $v \neq v_0$.
Hence we also have
\begin{equation} \label{eqn:BT'q}
 B(T', \textbf{\emph{d}}) =  \sum_{
        \scriptstyle p: V(T')\setminus \{v_0\} \longrightarrow \mathbb{Z}_{\geq 0}}
                 e_* [T_p , q_{\textbf{\emph{d}}}].
\end{equation}
In practice, this equation may be easier to use than equation \eqref{eqn:BT'}.

We first prove a result which is stronger than Proposition \ref{prop: cor 2} in the case $l=1$.
Let $\pi : { \overline {\cal M} }_{g, n+m+1} \to { \overline {\cal M} }_{g, n+m}$ be the map forgetting marked points labelled by the last frozen leg, i.e. $\pi$ is the map $\pi^{m,n}_1$ in Proposition \ref{prop: cor 2}.
Let $\tilde{\pi} : {\mathcal{A'}}_{g,n,m+1} \to {\mathcal{A'}}_{g,n,m}$ be the map forgetting the last frozen leg. Then we have
\begin{prop}\label{lem}
For any $m \geq 2 , n \geq 1 , g \geq 0$, and $T' \in {\mathcal{A'}}_{g,n,m+1}$, we have
\begin{equation}\label{push forward of B 1}
\pi_* B(T', \textbf{d})=  \sum_{\textbf{k} \, \in D_1(\textbf{d}) } B(\tilde{\pi}(T'), \textbf{k})
\end{equation}
for any $\vd \in (\mathbb{Z}_{\geq0})^n$.
\end{prop}

For any $T' \in {\mathcal{A}'}_{g,n,m}$, a vertex $v \in V(T')$ is called a {\it branching vertex} if there are at least two distinct positively directed half edges attached to $v$. If $n=1$, then $T'$ can have only one top vertex since every top vertex must have at least
one regular leg. In this case, $T'$  has the following shape
\begin{center}
\includegraphics[scale=0.6]{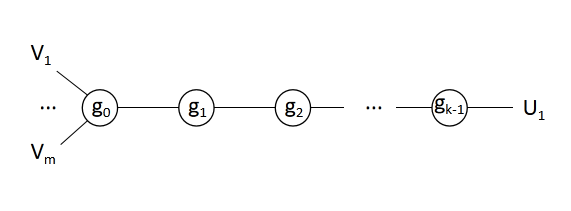}
\end{center}
where $U_1$ is the unique regular leg and $V_1, \ldots, V_m$ are frozen legs. In this case $T'$ has no branching vertex.
Note that a top vertex with at least two distinct regular legs
is considered as a branching vertex.  Hence if $n \geq 2$,  the set of branching vertices is non-empty and
 there exists a unique branching vertex, called the {\it first branching vertex},  whose level is less than the levels of all other branching vertices. The level of the first branching vertex is called the {\it branching height} of $T'$, which is denoted by $b(T')$.
 In particular, $b(T')=1$ if and only if the root of $T'$ is a branching vertex.

The following lemma allows us to prove Proposition \ref{lem} by induction on the branching height.
\begin{lem}\label{lem2}
Given any $N \geq 2$, if we have proved Proposition \ref{lem} for $n=N$ and $b(T')=1$, then Proposition \ref{lem} is also true for $n=N$ and $b(T') > 1$.
\end{lem}

\noindent
{\bf Proof}:
We prove this lemma by induction on $b(T')$. Assume Proposition \ref{lem} is true for $n=N$ and $b(T')=a \geq 1$.
We want to show that  Proposition \ref{lem} is also true for $n=N$ and $b(T')=a+1$.

Let $T' \in \mathcal{A}'_{g,n,m+1}$ with $b(T')=a+1$. The root $v_0$ of $T'$ can not be a branching vertex since $b(T')>1$. Hence
there is a unique vertex $v_1$ which is connected to $v_0$ by a single edge. Then $T'$ has the following shape
\begin{center}
\includegraphics[scale=0.6]{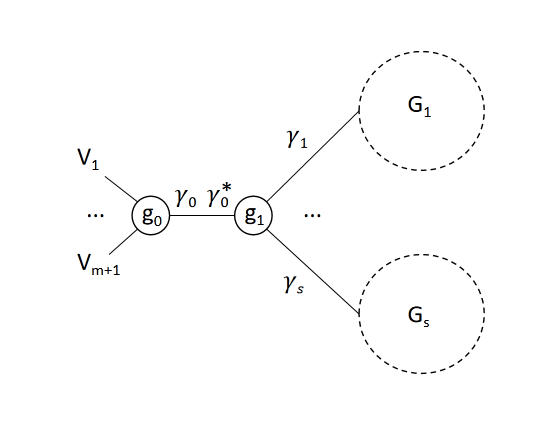}
\end{center}
 where $g_0$ and $g_1$ are the genera of $v_0$ and $v_1$ respectively, $\{ V_1, \ldots, V_{m+1}\}$ is the set of all frozen legs,
$\gamma_0$ is the unique positively directed half edge on $v_0$,
$\{ \gamma_1,  \ldots, \gamma_s\}$ is the set of all positively directed half edges on $v_1$,
$G_i$ is the part of the tree connected to $\gamma_i$ for $1 \leq i \leq s$. Note that $G_i$ may be empty, in which case $\gamma_i$ is just a regular leg.


Note that any $T \in {\mathcal{A}}(T')$ is obtained from $T'$ by adding certain numbers of extra legs to vertices of
$T'$ which are not the root. Let $\{ W_1, \ldots, W_{q_0+1}\}$ be the set of extra legs on $T$ attached to the vertex $v_1$.
Note that for each
$1 \leq i \leq s$, a vertex in $G_i$ is also a vertex in $T$.
Let $\widetilde{G}_i$ be the tree obtained from $G_i$ by adding same sets of extra legs to vertices as in corresponding vertices
on $T$.
Given any $\vd \in (\mathbb{Z}_{\geq 0})^n$, it induces a $\psi$-decoration
$q_{\vd}: H(T) \longrightarrow \mathbb{Z}_{\geq0}$ as given by equation \eqref{eqn:qd}. Note that $H(\widetilde{G}_i) \subset H(T)$. So
the restriction of $q_{\vd}$ to $H(\widetilde{G}_i)$ gives a $\psi$-decoration on the tree $\widetilde{G}_i$. Let $[G_i]_{T, \vd}$ be the
tautological class corresponding to this $\psi$-decorated graph. Let $q_j = q_{\vd}(\gamma_j)$ for $j=1, \ldots,s$.
Recall that the map $\pi$ forgets the marked point labeled by the frozen leg $V_{m+1}$. We have
\begin{eqnarray}
&& \pi_* (B(T', \textbf{\emph{d}})) \nonumber \\
&=&
 \sum_{\scriptstyle T \in {\mathcal{A}}(T')}
 \pi_* e_*\bigg( <\Psi^{q_0} (\gamma_0) \prod_{i=1}^{m+1} V_i  >_{g_0} <\gamma^*_0 \prod_{i=1}^{q_0+1} W_i  \prod_{j=1}^s \Psi^{q_j} (\gamma_j)>_{g_1} \prod_{i=1}^s [G_i]_{T, \vd}\bigg) \nonumber \\
&=&
 \sum_{\scriptstyle T \in {\mathcal{A}}(T')} <\Psi^{q_0 -1} (\gamma_0) \prod_{i=1}^{m} V_i  >_{g_0}
  e_*\bigg(  <\gamma^*_0 \prod_{i=1}^{q_0+1} W_i  \prod_{j=1}^s \Psi^{q_j} (\gamma_j)>_{g_1} \prod_{i=1}^s [G_i]_{T, \vd}\bigg),
  \label{eqn:BTGi}
\end{eqnarray}
where in the second equality, we have used string equation to get rid of $\pi_*$ and also used equation \eqref{eqn:e*T} to
compute $e_*$. In case $G_i$ is empty, $[G_i]_{T, \vd}$ should be omitted in the above equation.

Let $\tilde{e}$ be the map forgetting marked points labeled by extra legs on $\widetilde{G}_1, \ldots, \widetilde{G}_s$.
Let $\pi_1$ be the map forgetting marked points labeled by $W_1, \ldots, W_{q_0}$,
and $\pi_2$ the map forgetting marked points labeled by $W_{q_0+1}$. Then
$e=\pi_1 \pi_2 \tilde{e}$.

Let $T''$ be the rooted tree obtained from $T$ by first removing $v_0$, $\gamma_0$, and all extra legs attached to vertices
which are not $v_1$, then set $v_1$ to be the root and
regard $\gamma_0^*, W_1, \ldots, W_{q_0+1}$ as frozen legs of $T''$. Note that there is no regular legs attached to $v_0$ since
$v_0$ is not a branching vertex. Hence all regular legs of $T$ are attached to vertices in $T''$.
Then equation \eqref{eqn:BTGi} becomes
\begin{eqnarray}
\pi_* (B(T', \textbf{\emph{d}}))
&=&
 \sum_{q_0 \geq 1} <\Psi^{q_0 -1} (\gamma_0) \prod_{i=1}^{m} V_i  >_{g_0} (\pi_1)_* (\pi_2)_*
 \bigg( \sum_{\scriptstyle \widetilde{T} \in {\mathcal{A}}(T'')} (e_{\widetilde{T}})_* [\widetilde{T}, q_{\vd}] \bigg),
  \label{eqn:BTtilde}
\end{eqnarray}
where $e_{\widetilde{T}}$ is the map forgetting all marked points labelled by extra legs on $\widetilde{T}$.
For $\widetilde{T} \in {\mathcal{A}}(T'')$, $e_{\widetilde{T}}$ is the same map as $\tilde{e}$ above.

Since $b(T'')=a$, by induction hypothesis,
\[ (\pi_2)_*
 \bigg( \sum_{\scriptstyle \widetilde{T} \in {\mathcal{A}}(T'')} (e_{\widetilde{T}})_* [\widetilde{T}, q_{\vd}] \bigg)
 = (\pi_2)_* B(T'', \vd) = \sum_{\vk \in D_1(\vd)} B(\tilde{\pi}_2 (T''), \vk),
 \]
where $\tilde{\pi}_2 (T'')$ is the tree obtained from $T''$ by forgetting the frozen leg
$W_{q_0+1}$.
Hence equation \eqref{eqn:BTtilde} becomes
\begin{eqnarray*}
\pi_* (B(T', \textbf{\emph{d}}))
&=& \sum_{\vk \in D_1(\vd)} \,\,
 \sum_{q_0 \geq 1} <\Psi^{q_0 -1} (\gamma_0) \prod_{i=1}^{m} V_i  >_{g_0} (\pi_1)_* B(\tilde{\pi}_2 (T''), \vk) \nonumber \\
&=& \sum_{\vk \in D_1(\vd)} \,\,
 \sum_{q_0 \geq 1} <\Psi^{q_0 -1} (\gamma_0) \prod_{i=1}^{m} V_i  >_{g_0} (\pi_1)_*
        \bigg( \sum_{\widetilde{T} \in \mathcal{A} (\tilde{\pi}_2 (T''))}
              (e_{\widetilde{T}})_* [\widetilde{T}, q_{\vk}] \bigg)\nonumber \\
&=& \sum_{\vk \in D_1(\vd)} \,\, \sum_{T \in \mathcal{A} (\tilde{\pi} (T'))}
               (e_T)_* [T, q_{\vk}] \nonumber \\
&=& \sum_{\vk \in D_1(\vd)} B(\tilde{\pi} (T'), \vk) .
\end{eqnarray*}
This finishes the proof of the lemma.
$\Box$

We are now ready to prove Proposition \ref{lem}.

\vspace{10pt}
\noindent
{\bf Proof of Proposition \ref{lem}}:
We prove this Proposition by induction on $n$.

Assume $n=1$. Then $\textbf{\emph{d}}$ has only one component and we simply write $\textbf{\emph{d}}$ as $d \in \mathbb{Z}_{\geq 0}$.
Any $T' \in \mathcal{A}'_{g,1,m+1}$ can have only one top vertex, and it
must have the following shape
\begin{center}
\includegraphics[scale=0.6]{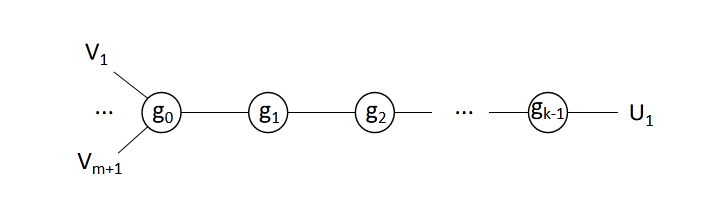}
\end{center}
 where $k$ is the number of vertices on $T'$, $g_i$ is the genus for the vertex $v_i$ for $0 \leq i \leq k-1$,
$\{V_1, \ldots, V_{m+1}\}$ is the set of frozen legs on $T'$, $U_1$ is the only regular leg.
Let $\gamma_i$ be the unique positively directed half edge on the vertex $v_i$ for $0 \leq i \leq k-1$.
Note that $\gamma_{k-1}=U_1$.
Any $T \in \mathcal{A}(T')$ is uniquely determined by $k-1$ nonnegative integers $q_i$ such that $q_i +1$ is the
number of extra legs in $T$ attached to the vertex $v_i$ for $i=1, \ldots, k-1$.
Hence by equation \eqref{eqn:BT'q}, we have
\begin{eqnarray}
&& \pi_* (B(T', d)) \nonumber \\
&=& \sum_{q_1, \ldots, q_{k-1} \in \mathbb{Z}_{\geq 0}}
            \pi_*  e_* \bigg( <V_1 \ldots V_{m+1} \Psi^{q_1}(\gamma_0)>_{g_0}
                \prod_{i=1}^{k-1}  <\gamma^*_{i-1}  \Psi^{q_{1+i}}(\gamma_i) \prod_{j=1}^{q_i +1} W^i_j >_{g_i}
             \bigg)
        \nonumber
\end{eqnarray}
where we have set $q_k = d$, $\{ W^i_1, \ldots, W^i_{q_i +1} \}$ is the set of extra legs on the vertex $v_i$
for $1 \leq i \leq k-1$. Using the string equation to get rid of $\pi_*$ and $e_*$ on the right hand side of this equation,
we obtain
\begin{eqnarray}
&& \pi_* (B(T', d)) \nonumber \\
&=& \sum_{q_1, \ldots, q_{k-1} \in \mathbb{Z}_{\geq 0}}
             <V_1 \ldots V_{m} \Psi^{q_1-1}(\gamma_0)>_{g_0}
                \prod_{i=1}^{k-1}  <\gamma^*_{i-1}  \Psi^{q_{1+i}-q_i -1}(\gamma_i) >_{g_i}.
       \label{B m+1 1}
\end{eqnarray}
Note that a necessary condition for the summands on the right hand side of this equation
to be non-zero is that
\[ q_1 \geq 1, \hspace{10pt} {\rm and \,\,\,} q_{1+i} \geq 1+ q_i {\rm \,\,\, for \,\,\, } 1 \leq i \leq k-1.\]
On the other hand, by equation \eqref{eqn:BT'q},
\begin{eqnarray}
&&  B(\tilde{\pi}(T'), d-1) \nonumber \\
&=& \sum_{p_1, \ldots, p_{k-1} \in \mathbb{Z}_{\geq 0}}
    e_* \bigg( <V_1 \ldots V_{m} \Psi^{p_1}(\gamma_0)>_{g_0}
                \prod_{i=1}^{k-1}  <\gamma^*_{i-1}  \Psi^{p_{1+i}}(\gamma_i) \prod_{j=1}^{p_i +1} W^i_j >_{g_i}
             \bigg),
     \nonumber
\end{eqnarray}
where $\{W^i_1, \ldots, W^i_{1+p_i}\}$ in the set of extra legs on the vertex $v_i$ in $T \in \mathcal{A}(\tilde{\pi}(T'))$
for $1 \leq i \leq k-1$,
$p_k=d-1$. Using the string equation to get rid of $e_*$, we have
\begin{eqnarray}
&& B(\tilde{\pi}(T'), d-1) \nonumber \\
&=& \sum_{p_1, \ldots, p_{k-1} \in \mathbb{Z}_{\geq 0}}
             <V_1 \ldots V_{m} \Psi^{p_1}(\gamma_0)>_{g_0}
                \prod_{i=1}^{k-1}  <\gamma^*_{i-1}  \Psi^{p_{1+i}-p_i -1}(\gamma_i) >_{g_i}
                \label{B m 1}.
\end{eqnarray}
For the summands on the right hand side of this equation to be non-zero, $p_i$ must satisfy the condition
\[ p_1 \geq 0, \hspace{10pt} {\rm and \,\,\,} p_{1+i} \geq 1+ p_i {\rm \,\,\, for \,\,\, } 1 \leq i \leq k-1.\]
Setting
$p_i = q_i - 1$ for all $i=1, \cdots, k$, we see the right hand sides of equations \eqref{B m+1 1} and \eqref{B m 1}
are the same. Hence
\[ \pi_* (B(T', d)) = B(\tilde{\pi}(T'), d-1) .\]
This proves Proposition \ref{lem} for $n=1$.

Assume  Proposition \ref{lem} holds for all $n$ such that $1 \leq n \leq N$ where $N$ is an arbitrary positive integer.
We want to prove this proposition also holds
 for $n=N+1$.

We first observe that if $T'$ has only one vertex, then $\mathcal{A}(T')=\{T'\}$ and equation \eqref{push forward of B 1}
is just the string equation. So we may assume $T'$ has more than one vertices.

 By Lemma \ref{lem2}, we only need to show this proposition holds
 for all $T' \in \mathcal{A}'_{g, N+1, m+1}$ with branching height equal to 1, i.e.
 there are at least two positively directed half edges attached to the root
 $v_0$ of $T'$.
Such $T'$ has the following shape
\begin{center}
\includegraphics[scale=0.6]{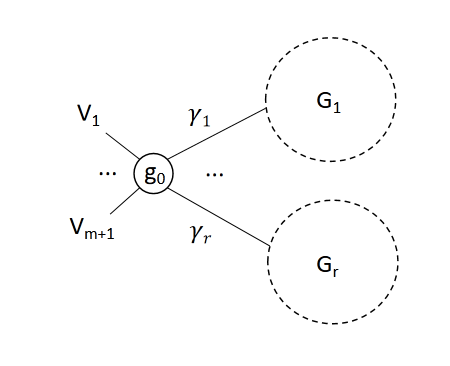}
\end{center}
where $g_0$ is the genus of the root $v_0$, $\{V_1, \ldots, V_{m+1}\}$ is the set of all frozen legs on $T'$,
$\{ \gamma_1, \ldots, \gamma_r \}$ is the set of all positively directed half edges on $v_0$ where $r \geq 2$,
and $G_i$ is the  part of tree connected to $\gamma_i$ for $1 \leq i \leq r$. Note that
 $G_i$ is an empty tree if $\gamma_i$ is a regular leg. Without loss of generality, we may assume
 that $G_1, \cdots, G_s$ are non-empty, and $\gamma_{s+1}, \cdots \gamma_r$ are regular legs of $T'$
 for some $1 \leq s  \leq r$. For $1 \leq i \leq s$,
 let $v_i$ be the unique vertex on $G_i$ which is connected to $v_0$ by a single edge.

Every $T \in {\mathcal{A}}(T')$ is uniquely determined by a function
\[ p: V(T') \setminus \{v_0\} \longrightarrow \mathbb{Z}_{\geq 0} \]
such that $p(v)+1$ is the number of extra legs on $T$ attached to vertex $v \neq v_0$.
Since $V(T') \setminus \{v_0\}$ is equal to the disjoint union of $V(G_i)$ for $i=1, \ldots, s$,
the map $p$ is uniquely determined by a set of $s$ functions
\[ p_i: V(G_i) \longrightarrow \mathbb{Z}_{\geq 0}, \,\,\, i=1, \ldots, s \]
such that $p_i$ is the restriction of $p$ to $V(G_i)$.
Let $\widetilde{G}_i$ be the tree obtained from $G_i$ by
adding $p_i(v)+1$ extra legs to the vertex $v \in V(G_i)$.
Let $I_i$ be the set of indices for regular legs which lie in $G_i$. Since  $r \geq 2$ and each $I_i$ is
nonempty for $1 \leq i \leq s$, the number of regular legs in each $G_i$ is less than or equal to $N$.

Given a weight vector $\vd \in (\mathbb{Z}_{\geq 0})^{N+1}$ for $T'$,
let $\vd_i$ be the weight vector on $G_i$ obtained by restricting $\vd$ to regular legs in $G_i$.
Together with $p_i$, this induces a $\psi$-decoration $q_{\vd_i}$ on $\widetilde{G}_i$ in the same way as in equation
\eqref{eqn:qd}.
The tautological class corresponding to this $\psi$-decorated graph is denoted by $[G_i]_{p_i, \vd_i}$.
Let $q_i = p_i(v_i)$ for $i=1, \ldots, s$ and $q_j$ equal to the weight of the regular leg $\gamma_j$ as given by the
weight vector $\vd$ for $j=s+1, \ldots, r$. Then by equation \eqref{eqn:BT'q}, we have
\begin{eqnarray}
\pi_* (B(T', \textbf{\emph{d}}))
&=& \sum_{\scriptstyle p_1, \ldots, p_s} \pi_* e_* \bigg(
        <\prod_{l=1}^{m+1} V_l \prod_{k=1}^r \Psi^{q_k}(\gamma_k) >_{g_0}  \prod_{j=1}^s [G_j]_{p_j, \vd_j} \bigg),
        \nonumber
\end{eqnarray}
where the summation $\sum_{\scriptstyle p_1, \ldots, p_s}$ is over the set all
maps $p_1, \ldots, p_s$.

Using string equation to get rid of $\pi_*$ and applying equation \eqref{eqn:e*T}, we have
\begin{eqnarray}
\pi_* (B(T', \textbf{\emph{d}}))
&=& \sum_{\scriptstyle p_1, \ldots, p_s}  \sum_{i=1}^r
        <\Psi^{q_i -1}(\gamma_i) \prod_{l=1}^{m} V_l  \prod_{k=1, \atop k \neq i}^r \Psi^{q_k}(\gamma_k) >_{g_0}
                             \prod_{j=1}^s e_* [G_j]_{p_j, \vd_j} \nonumber \\
&=& \sum_{i=1}^s \sum_{\scriptstyle p_1, \ldots, \widehat{p_i}, \ldots,  p_s} \Phi_i \prod_{j=1, \atop j \neq i}^s e_* [G_j]_{p_j, \vd_j}
    + \sum_{i=s+1}^r \sum_{\scriptstyle p_1, \ldots, p_s} \Phi_i \prod_{j=1}^s e_* [G_j]_{p_j, \vd_j},
    \hspace{10pt}
                \label{eqn:BT'bPhi}
\end{eqnarray}
where
\begin{equation} \label{eqn:Phi1}
\Phi_i := \sum_{p_i} <\Psi^{q_i -1}(\gamma_i) \prod_{l=1}^{m} V_l  \prod_{k=1, \atop k \neq i}^r \Psi^{q_k}(\gamma_k) >_{g_0}
                      e_* [G_i]_{p_i, \vd_i}
\end{equation}
for $1 \leq i \leq s$, and
\begin{equation} \label{eqn:Phi2}
\Phi_i :=  <\Psi^{q_i -1}(\gamma_i) \prod_{l=1}^{m} V_l  \prod_{k=1, \atop k \neq i}^r \Psi^{q_k}(\gamma_k) >_{g_0}
\end{equation}
for $s+1 \leq i \leq r$.

Assume $1 \leq i \leq s$.
Let $T_i'$ be the tree obtained from $T'$ by removing all $G_j$ if $j \neq i$ and take
$\{ \gamma_j \mid 1 \leq j \leq r, j \neq i\} \bigcup \{ V_1, \ldots, V_{m+1}\}$ as the set of all frozen legs on $T_i'$.
Then the number of regular legs on $T_i'$ is equal to
the number of regular legs on $G_i$, which is less than or equal to $N$. By induction hypothesis,
\begin{equation} \label{eqn:Indn}
\pi_* B(T_i', \vd_i) \\
= \sum_{\vk_i \in D_1(\vd_i)} B(\tilde{\pi}(T_i'), \vk_i),
\end{equation}
where $\pi$ forgets the marked points labelled by the frozen leg $V_{m+1}$ on $T_i'$.

Since
$V(T_i') \setminus \{v_0\} = V(G_i)$, $\mathcal{A}(T_i')$ is one to one correspondent to the set of all maps
$p_i: V(G_i) \longrightarrow \mathbb{Z}_{\geq 0}$. Hence by equation \eqref{eqn:BT'q}, we have
\begin{eqnarray*}
\pi_* B(T_i', \vd_i) &=& \pi_* \sum_{p_i} e_* \bigg(
            <\Psi^{q_i}(\gamma_i) \prod_{l=1}^{m+1} V_l  \prod_{k=1, \atop k \neq i}^r \gamma_k >_{g_0}
                       [G_i]_{p_i, \vd_i} \bigg) \\
&=& \sum_{p_i} <\Psi^{q_i -1}(\gamma_i) \prod_{l=1}^{m} V_l  \prod_{k=1, \atop k \neq i}^r \gamma_k >_{g_0}
                      e_* [G_i]_{p_i, \vd_i},
\end{eqnarray*}
and
\begin{eqnarray*}
\sum_{\vk_i \in D_1(\vd_i)} B(\tilde{\pi}(T_i'), \vk_i)
 &=& \sum_{\vk_i \in D_1(\vd_i)}
   \sum_{p_i} e_*  \bigg( <\Psi^{q_i}(\gamma_i) \prod_{l=1}^{m} V_l  \prod_{k=1, \atop k \neq i}^r \gamma_k >_{g_0}
                      [G_i]_{p_i, \vk_i} \bigg).
\end{eqnarray*}
Hence equation \eqref{eqn:Indn} is equivalent to
\begin{eqnarray*}
&& \sum_{p_i} <\Psi^{q_i -1}(\gamma_i) \prod_{l=1}^{m} V_l  \prod_{k=1, \atop k \neq i}^r \gamma_k >_{g_0}
                      e_* [G_i]_{p_i, \vd_i} \\
 &=& \sum_{\vk_i \in D_1(\vd_i)}
   \sum_{p_i}  <\Psi^{q_i}(\gamma_i) \prod_{l=1}^{m} V_l  \prod_{k=1, \atop k \neq i}^r \gamma_k >_{g_0}
                      e_* [G_i]_{p_i, \vk_i}.
\end{eqnarray*}
 Note that $\gamma_k$ are legs on $T_i'$ for $k \neq i$. Let $\psi_{\gamma_k}$ be the corresponding $\psi$-classes.
 Multiplying both sides of the above equation by $$\prod_{k=1,\atop k \neq i}^r (\psi_{\gamma_k})^{q_k}$$ and using equation
 \eqref{11111}, we obtain
 \begin{eqnarray*}
\Phi_i
&=& \sum_{p_i} <\Psi^{q_i -1}(\gamma_i) \prod_{l=1}^{m} V_l  \prod_{k=1, \atop k \neq i}^r \Psi^{q_k}(\gamma_k) >_{g_0}
                      e_* [G_i]_{p_i, \vd_i} \\
 &=& \sum_{\vk_i \in D_1(\vd_i)}
   \sum_{p_i}  <\Psi^{q_i}(\gamma_i) \prod_{l=1}^{m} V_l  \prod_{k=1, \atop k \neq i}^r \Psi^{q_k}(\gamma_k) >_{g_0}
                      e_* [G_i]_{p_i, \vk_i},
\end{eqnarray*}
where $\Phi_i$ is defined by equation \eqref{eqn:Phi1} for $1 \leq i \leq s$.
Plugging this formula into equation \eqref{eqn:BT'bPhi}, we have
\begin{eqnarray*}
\pi_* (B(T', \textbf{\emph{d}}))
&=& \sum_{i=1}^s \sum_{\vk_i \in D_1(\vd_i)} \,\, \sum_{\scriptstyle p_1, \ldots,  p_s}
        < \prod_{l=1}^{m} V_l  \prod_{k=1}^r \Psi^{q_k}(\gamma_k) >_{g_0}
                      e_* [G_i]_{p_i, \vk_i} \prod_{j=1, \atop j \neq i}^s e_* [G_j]_{p_j, \vd_j} \\
&&    + \sum_{i=s+1}^r \sum_{\scriptstyle p_1, \ldots, p_s}
    <\Psi^{q_i-1}(\gamma_i) \prod_{l=1}^{m} V_l  \prod_{k=1, \atop k \neq i}^r \Psi^{q_k}(\gamma_k) >_{g_0} \prod_{j=1}^s e_* [G_j]_{p_j, \vd_j}  \\
&=& \sum_{\vk \in D_{1}(\vd)} B(\tilde{\pi}(T'), \vk).
\end{eqnarray*}
Hence equation \eqref{push forward of B 1} holds for $T'$.
This finishes the proof of the proposition.
$\Box$

An immediate consequence of Proposition \ref{lem} and equation \eqref{eqn:B=BT'} is the following
\begin{cor}\label{cor 1}
Under the same conditions as in Proposition \ref{lem}, we have
\[
\pi_* (B^{m+1}_{g, \textbf{d}})=  \sum_{\textbf{k} \in D_1(\textbf{d})}
     B^m_{g, \textbf{k}} \,.
\]
\end{cor}

We are now ready to prove the main result of this section.

\vspace{6pt}
\noindent
{\bf Proof of Proposition \ref{prop: cor 2}}:
We prove this proposition by induction on $l$. If $l=1$, this proposition is just Corollary \ref{cor 1}.
Assume the proposition holds for an arbitrary positive integer $l$.
We want to show it also holds when $l$ is replaced by $l+1$.

Since $\pi_{l+1}^{m,n} = \pi_{1}^{m,n} \pi_{l}^{m+1,n}$, by induction hypothesis and Corollary \ref{cor 1}, we have
\begin{eqnarray*}
(\pi_{l+1}^{m,n})_* (B^{m+l+1}_{g, \textbf{\emph{d}}})
&=& (\pi_{1}^{m,n})_* (\pi_{l}^{m+1,n})_* \,\, B^{m+l+1}_{g, \textbf{\emph{d}}}
= (\pi_{1}^{m,n})_* \bigg( \sum_{\vk \in D_l(\vd) } \frac{l!}{\prod^n_{i=1} (d_i - k_i)!} B^{m+1}_{g, \textbf{\emph{k}}} \bigg) \\
&=& \sum_{\vk \in D_l(\vd) }  \frac{l!}{\prod^n_{i=1} (d_i - k_i)!} \sum_{\textbf{\emph{h}} \in D_1(\vk)}
                 B^{m}_{g, \textbf{\emph{h}}}
= \sum_{\textbf{\emph{h}} \in D_{1+1}(\vd)} \frac{(l+1)!}{\prod^n_{i=1} (d_i - h_i)!} B^{m}_{g, \textbf{\emph{h}}},
\end{eqnarray*}
where the last equality follows from the following combinatorial identity
\[ \sum^n_{j=1} \frac{l!}{\prod^n_{i=1} (d_i - h_i - \delta_{ij})!}
=\frac{l!}{\prod^n_{i=1} (d_i - h_i)!} \sum_{j=1}^n (d_j-h_j)
= \frac{(l+1)!}{\prod^n_{i=1} (d_i - h_i)!}. \]
Here $\delta_{ij} = 1$ if $i=j$ and equal to $0$ otherwise. This completes the proof of the proposition.
$\Box$

\section{Proof of Theorem \ref{thm:m=2}}
\label{section induction of thm}

In this section, we will prove Theorem \ref{thm:m=2} by induction on $m$.
Lemma 7.2 in \cite{BS} says that we can reduce the Buryak-Shadrin conjecture to the case where all components
of $\vd$ are positive. Lemma 7.3 in \cite{BS} says we can reduce this conjecture to the case
where $|\vd|=2g+m-1$. Therefore to prove Theorem \ref{thm:m=2}, we only need to show that
this conjecture can be reduced to the case where $m=2$. So we only need to prove the following
\begin{prop}\label{induction step}
Given $m \geq 2$ and $g \geq 0$, if $B^{m'}_{g', \vd} =0$ for all $2 \leq m' \leq m$, $g' \leq g$, $n \geq 1$,
and $\vd \in (\mathbb{Z}_{\geq 0})^n$ with $|\vd| \geq 2g'+m'-1$, then we have
$B^{m+1}_{g, \vd} =0$ for all $n \geq 1$ and $\vd \in (\mathbb{Z}_{\geq 0})^n$ with $|\vd| \geq 2g+m$.
\end{prop}

\noindent
{\bf Proof}:
The proposition follows if we can prove that for any fixed $g \geq 0$, $n \geq 1$, $\vd \in (\mathbb{Z}_{\geq 0})^n$
with $|\textbf{\emph{d}}| \geq 2g+m$,
$ B^{m+1}_{g, \textbf{\emph{d}}} = B^{m}_{g, (\textbf{\emph{d}}, 0)}$ under the assumption of the proposition.
This is equivalent to show that  $B^{m}_{g, (\textbf{\emph{d}}, 0)} - B^{m+1}_{g, \textbf{\emph{d}}} =0$. For this purpose, we first compute
\begin{equation*} 
B^{m}_{g, (\textbf{\emph{d}}, 0)} = \sum_{T \in {\mathcal{A}}_{g,n+1,m}}  (-1)^{|E(T)|} e_* [T , q_{(\textbf{\emph{d}}, 0)}].
\end{equation*}

Let $\{ U_1, \ldots, U_{n+1}\}$ be the set of all regular legs on a tree $T \in {\mathcal{A}}_{g,n+1,m}$. The weights of
the first $n$ regular legs are given by $\vd$ and
 the weight of the last regular leg $U_{n+1}$ is $0$.
Assume $U_{n+1}$ is attached to a vertex $v_1$. If $v_1$ is not
the root $v_0$, there exists a unique vertex $v_2$ which is connected to $v_1$ by a single edge $e_{1,2}$ such that the level of $v_2$ is strictly
less than the level of $v_1$. Cutting $T$ along the edge $e_{1,2}$, we obtain two non-empty connected subtrees $T_1$ and $T_2$
such that $v_i$ lies on $T_i$ for $i=1, 2$.
Set
\[ I := \{ i \mid 1 \leq i \leq n, U_i \in L(T_1) \}, \hspace{20pt} J := \{ j \mid U_j \in L(T_2) \}. \]
Let $\vd_I$ be the weight vector whose components are weights of regular legs on $T_1$ and
 $\vd_J$  the weight vector whose components are weights of regular legs on $T_2$.
Let $g(T_i)$ be the genus of $T_i$. Then $g(T_1)+g(T_2)=g$.
Let $p$ be the non-negative integer such that $p+1$ is the number of extra legs of $T$ attached to the vertex $v_1$.
Since $|\textbf{\emph{d}}_I|+ |\textbf{\emph{d}}_J| = |\vd|$,
 one and only one of the following two inequalities must hold:
\begin{equation} \label{eqn:d>g2}
 |\textbf{\emph{d}}_I| \geq 2g(T_1) + 2 +p
 \end{equation}
 or
 \begin{equation} \label{eqn:d>g1}
|\textbf{\emph{d}}_J| + p \geq |\textbf{\emph{d}}| - 2g(T_1) - 1.
\end{equation}

Let $S_1$ be the set of all $T \in {\mathcal{A}}_{g,n+1,m}$ such that $v_1 \neq v_0$ and inequality \eqref{eqn:d>g2} holds,
and  $S_2$ the set of all $T \in {\mathcal{A}}_{g,n+1,m}$ such that $v_1 \neq v_0$ and inequality \eqref{eqn:d>g1} holds.
Moreover, let $S_0$ be the set of all $T \in {\mathcal{A}}_{g,n+1,m}$ such that  $v_1 = v_0$.
Then we have a decomposition of ${\mathcal{A}}_{g,n+1,m}$ as a disjoint union of three subsets
\[ {\mathcal{A}}_{g,n+1,m} = S_0 \coprod S_1 \coprod S_2. \]
Let
\begin{equation} \label{eqn:Theta}
\Theta_i := \sum_{T \in S_i}  (-1)^{|E(T)|} e_* [T , q_{(\textbf{\emph{d}}, 0)}]
\end{equation}
for $i=0, 1, 2$. Then we have
\begin{equation} \label{eqn:BTheta}
B^{m}_{g, (\textbf{\emph{d}}, 0)} = \Theta_0 + \Theta_1 + \Theta_2 .
\end{equation}

For $T \in S_0$, the last regular leg $U_{n+1}$ is attached to the root and its weight is $0$. So in the definition
of tautological class $e_* [T , q_{(\textbf{\emph{d}}, 0)}]$, $U_{n+1}$ behaves in the same way as a frozen leg.
Hence we have
\begin{equation}\label{hhh 1}
\Theta_0 = \sum_{T \in {\mathcal{A}}_{g,n,m+1}}  (-1)^{|E(T)|} e_* [T , q_{\textbf{\emph{d}}}] = B^{m+1}_{g, \textbf{\emph{d}}}.
\end{equation}
Therefore we only need to prove $\Theta_1 = \Theta_2 =0$ under the assumption of the proposition.

For $T \in S_1 \bigcup S_2$, it can be represented as two non-empty subtrees $T_1$ and $T_2$ connected by an edge $e_{1,2}$
as described above.  This decomposition of $T$ is illustrated by the following picture

\vspace{40pt}
\begin{center}
\includegraphics[scale=0.6]{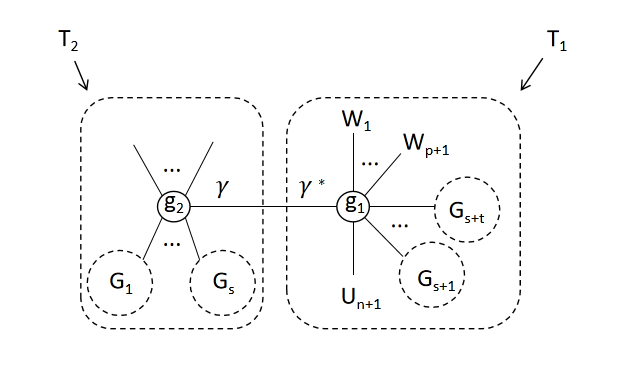},
\end{center}
 where $g_i$ is the genus of the vertex $v_i$ for $i=1, 2$,  $\gamma$ and $\gamma^*$ are respectively positively directed and negatively directed half edges of $e_{1,2}$, $W_1, \ldots, W_{p+1}$ are extra legs on the vertex $v_1$,
 $G_1, \ldots, G_{s+t}$ are possibly subtrees of $T$ connected to vertices $v_2$ and $v_1$ by some half edges.
 Note that the root $v_0$ lies in $T_2$.

Let $q(T_i)$ be the restriction of the $\psi$-decoration $q_{(\textbf{\emph{d}}, 0)}$  to $T_i$
for $i=1, 2$.
Let $[T_i, q(T_i)]$ for $i=1, 2$ be the corresponding tautological classes written in the form of equation \eqref{eqn:taut}.
Then the representation of $[T, q_{(\textbf{\emph{d}}, 0)}]$ in the form of equation \eqref{eqn:taut}
can be written as
\[ [T, q_{(\textbf{\emph{d}}, 0)}] = [T_1, q(T_1)] [T_2, q(T_2)].\]
Note that the right hand side of this equation  is not a product of two tautological classes.
It should be interpreted as a formal product as given in equation \eqref{eqn:taut}.
More precisely, the formal product
$[T_1, q(T_1)] [T_2, q(T_2)]$ means the push forward of the
class
\[ [T_1, q(T_1)] \otimes [T_2, q(T_2)]
\in H^*(\overline{\mathcal M}_{g(T_1), \, |L(T_1)|}) \otimes H^*(\overline{\mathcal M}_{g(T_2), \, |L(T_2)|})\]
under the gluing map
\[ \overline{\mathcal M}_{g(T_1), \, |L(T_1)|} \times \overline{\mathcal M}_{g(T_2), \, |L(T_2)|}
  \longrightarrow  \overline{\mathcal M}_{g(T_1)+g(T_2), \, |L(T_1)| + |L(T_2)| -2} \]
which glues the marked points labelled by $\gamma*$ and $\gamma$ together.
Hence for $i=1, 2$, we have
\begin{eqnarray*}
 \Theta_i
 &=& \sum_{T \in S_i} (-1)^{|E(T_1)|+|E(T_2)|+1} e_* \left( [T_1, q(T_1)] [T_2, q(T_2)] \right)  \\
 &=& \sum_{T \in S_i} (-1)^{|E(T_1)|+|E(T_2)|+1}   e_* [T_1, q(T_1)] \,\, e_*[T_2, q(T_2)].
\end{eqnarray*}

Note that $T_2$ is naturally a rooted tree with root $v_0$. We can consider $\gamma$ as a new regular leg on $T_2$ with weight equal to $p$ and treat it as the last regular leg on $T_2$.
We can construct a rooted tree $\tilde{T}_1$  from $T_1$ by setting
$v_1$ as the root of $\tilde{T}_1$ and consider $\gamma^*$, $U_{n+1}$ and all extra legs of $T$ attached to $v_1$ as frozen legs on $\tilde{T}_1$, where $\gamma^*$ and $U_{n+1}$ are treated as the first two frozen legs.
Let $n_i$ be the number of regular legs of $T$ which lie in $T_i$ for $i=1,2$.
Then
\[ \tilde{T}_1 \in \mathcal{A}_{g(T_1),n_1-1,p+3}, \hspace{10pt}   T_2  \in \mathcal{A}_{g(T_2),n_2+1,m}.\]
Let $\pi_{p+1}$ be the forgetful map which forgets marked points labelled by the last $p+1$ frozen legs on $\tilde{T}_1$.
Then
\[ \Theta_i  = \sum_{T \in S_i} (-1)^{|E(T_1)|+|E(T_2)|+1} (\pi_{p+1})_* e_* [\tilde{T}_1, q(T_1)] \,\, e_*[T_2, q(T_2)].\]

If $T \in S_1$, then total weights on $\tilde{T}_1$ is $|\vd_I| \geq 2 g(T_1)+p+2$. We have
\begin{eqnarray*}
\Theta_1 &=&  \sum_{p \geq 0} \sum_{0 \leq g' \leq g, \atop 0 \leq n' \leq n}
        \sum_{J \subset [n],  |J|=n', \atop |\vd|-|\vd_J| \geq 2 (g-g') +p+2}
        \sum_{T_2 \in \mathcal{A}_{g', n'+1,m}} (-1)^{|E(T_2)|+1} e_*[T_2, q_{(\vd_{J}, p)}]   \\
&& \hspace{100pt}
   \sum_{\tilde{T}_1 \in \mathcal{A}_{g-g',n-n',p+3}} (-1)^{|E(\tilde{T}_1)|} (\pi_{p+1})_* e_* [\tilde{T}_1, q_{\vd_I}],
\end{eqnarray*}
where $[n]:= \{1, 2, \ldots n\}$, $I=[n] \setminus J$, $\vd_I$ and $\vd_J$ mean the restrictions of $\vd$ to indices $I$ and $J$
respectively. 
If $n'=n$, $I$ is an empty set and $(\pi_{p+1})_* e_* [\tilde{T}_1, q_{\vd_I}]=0$ by the string equation.
If $n' < n$, the summation in the second line of the above formula is
\begin{eqnarray*}
\sum_{\tilde{T}_1 \in \mathcal{A}_{g-g',n-n',p+3}} (-1)^{|E(\tilde{T}_1)|} (\pi_{p+1})_* e_* [\tilde{T}_1, q_{\vd_I}]
&=& (\pi_{p+1})_* B^{p+3}_{g-g', \vd_I}.
\end{eqnarray*}
By Proposition \ref{prop: cor 2}, the right hand side of this equation is equal to a linear combination of
$B^{2}_{g-g', \vk}$ with $|\vk|=|\vd_{I}|-p-1 \geq 2(g-g')+1$. Hence it is equal to $0$ by assumption of Proposition \ref{induction step}.
This proves $\Theta_1=0$.

If $T \in S_2$, then total weights on $T_2$ is $|\vd_J| + p \geq |\vd|- 2 g(T_1)-1$. We have
\begin{eqnarray*}
\Theta_2 &=&  \sum_{p \geq 0} \sum_{0 \leq g' \leq g, \atop 0 \leq n' \leq n}
        \sum_{I \subset [n],  |I|=n', \atop p-|\vd_I| \geq -2g'-1}
        \sum_{\tilde{T}_1 \in \mathcal{A}_{g', n', p+3}} (-1)^{|E(\tilde{T}_1)|+1} (\pi_{p+1})_* e_*[\tilde{T}_1, q_{\vd_{I}}]   \\
&& \hspace{100pt}
   \sum_{T_2 \in \mathcal{A}_{g-g',n-n'+1,m}} (-1)^{|E(T_2)|}  e_* [T_2, q_{(\vd_J,p)}],
\end{eqnarray*}
where $J=[n] \setminus I$. The summation in the second line of the above formula is
\begin{eqnarray*}
\sum_{T_2 \in \mathcal{A}_{g-g',n-n'+1,m}} (-1)^{|E(T_2)|}  e_* [T_2, q_{(\vd_J,p)}]
&=&  B^{m}_{g-g', (\vd_J, p)}=0,
\end{eqnarray*}
where the last equality follows from the assumption of Proposition \ref{induction step} and the fact
\[ |\vd_J|+p = |\vd| - |\vd_I| +p \geq (2g+m) - 2 g' -1 = 2(g-g') +m-1. \]
This shows $\Theta_2=0$.
The proof of the proposition is thus finished.
$\Box$

This completes the proof of Theorem \ref{thm:m=2}.

\section{Proof of Theorem \ref{main theorem}}
\label{section computation of base case}

By equation \eqref{eqn:Bmgd},  $B^m_{g, \vd}$ is given by a summation of $(-1)^{|E(T)|} e_* [T, q_{\vd}]$ for $T \in \mathcal{A}_{g,n,m}$.
There are many redundant terms in this summation which are obviously equal to $0$ by the string equation. More precisely,
for any vertex $v \in V(T)$, let $p(v)$ be the number of extra legs attached to $v$, then
by equations \eqref{string equation} and \eqref{eqn:e*T},
\[ p(v) \leq \sum_{h \in H_v} q_{\vd}(h) \]
if $e_* [T, q_{\vd}] \neq 0$. On the other hand the total power of $\psi$-classes on a vertex $v$ should not exceed dimension of the moduli space
$\dim_{\mathbb{C}} \overline{\mathcal{M}}_{g_v, |H_v|} = 3 g_v - 3 + |H_v|$. Hence we have
\[ \sum_{h \in H_v} q_{\vd}(h)  \leq 3 g_v - 3 + |H_v \setminus L^0(T) | + p(v) \]
if $[T, q_{\vd}] \neq 0$.
Therefore to compute $B^m_{g, \vd}$, we only need to consider $T$ which satisfies the inequality
\begin{equation} \label{eqn:pvrange}
\sum_{h \in H_v} q_{\vd}(h) - 3 g_v + 3 - |H_v \setminus L^0(T)| \leq p(v) \leq \sum_{h \in H_v} q_{\vd}(h)
\end{equation}
for all $v \in V(T)$. This inequality puts strong restrictions for the number of extra legs on each vertex and
therefore reduces the number of trees $T$ to be considered in computing $B^m_{g, \vd}$.
To find non-trivial contributions to $B^m_{g, \vd}$, it is easier to first consider top vertices and the root
where it is easier to determine the number of extra legs.

By Theorem \ref{thm:m=2}, to prove Theorem \ref{main theorem}, we only need to check
the case $m=2$ and $|\vd|=3$ with all components positive since the genus $g=1$. Since the $n=1$ case has been proved in \cite{BS},
we only need to show
$B^2_{1, (2,1)}=0$ and $B^2_{1, (1,1,1)}=0$.
To prove these two equalities, we will use the following basic facts: $\psi_i=0$ on $\overline{\cal{M}}_{0,3}$ since
$\dim \overline{\cal{M}}_{0,3}=0$, and $\psi_1$ on $\overline{\cal{M}}_{1,1}$ is given by
\[  <\Psi(x_1)>_1=\frac{1}{12}<x_1 \, \gamma \, \gamma^*>_0, \]
 where $x_1$ represents the leg corresponding to the only marked point
 on each curve in $\overline{\cal{M}}_{1,1}$ (see, for example, page 385 in \cite{ACG}).
Repeatedly using the pullback
formula for $\psi$-classes under the forgetful map $\overline{\cal{M}}_{g,n+1} \longrightarrow \overline{\cal{M}}_{g,n}$
(see, for example,  Lemma 4.28 on page 581 in \cite{ACG}), it is straightforward
to check the following equations:
\begin{eqnarray}
&& <\Psi(x_1) x_2, \ldots, x_{a}>_0 \nonumber \\
&=& \sum_{k=1}^{a-3}  \sum_{2 \leq i_1 < \ldots < i_k \leq a-2} <x_1 x_{i_1} \ldots x_{i_k} \gamma>_0
         \, <\gamma^* x_2 \ldots \hat{x}_{i_1} \ldots \hat{x}_{i_k} \ldots x_{a}>_0
    \label{psi genus 0}
\end{eqnarray}
on $\overline{\cal{M}}_{0,a}$ with $a \geq 4$, and
\begin{eqnarray}
&& <\Psi(x_1) x_2, \ldots, x_{a}>_1   \nonumber \\
&=&  \sum_{k=1}^{a-1} \sum_{2 \leq i_1 < \ldots < i_k \leq a} <x_1 x_{i_1} \ldots x_{i_k} \gamma>_0
            \, <\gamma^* x_2 \ldots \hat{x}_{i_1} \ldots \hat{x}_{i_k} \ldots x_a>_1 \nonumber \\
&&  \hspace{80pt} + \frac{1}{12} <x_1 x_2 \ldots x_a \gamma \gamma^*>_0
    \label{psi genus 1}
\end{eqnarray}
on $\overline{\cal{M}}_{1,a}$ with $a \geq 1$.
Here $x_i$ for $i=1, \ldots, a$ represent legs corresponding to marked points on curves in $\overline{\cal{M}}_{g,a}$.
Multiplying both sides of equations \eqref{psi genus 0} and \eqref{psi genus 1} by more $\psi$-classes
using equation \eqref{11111} and applying equations \eqref{psi genus 0} and \eqref{psi genus 1} repeatedly,
we can obtain formulas which represent higher powers of $\psi$-classes on
genus-$0$ and genus-$1$ moduli spaces in terms of pure boundary classes which do not involve any $\psi$-classes.
Such formulas allow us to get rid of all $\psi$-classes in $B^2_{1, (2,1)}$ and $B^2_{1, (1,1,1)}$.
Then we can show $B^2_{1, (2,1)}$ and $B^2_{1, (1,1,1)}$
are $0$ by using
the well known genus-0 equation
\begin{eqnarray}
&& \sum_{k=0}^{a-4} \sum_{5 \leq i_1 < \ldots < i_k \leq a} <x_1 x_2 x_{i_1} \ldots x_{i_k} \gamma>_0 \, <\gamma^* x_3 x_4 x_5 \ldots \hat{x}_{i_1} \ldots \hat{x}_{i_k} \ldots x_a >_0 \nonumber \\
&=&
\sum_{k=0}^{a-4} \sum_{5 \leq i_1 < \ldots < i_k \leq a} <x_1 x_3 x_{i_1} \ldots x_{i_k} \gamma>_0 \, <\gamma^* x_2 x_4 x_5 \ldots \hat{x}_{i_1} \ldots \hat{x}_{i_k} \ldots x_a >_0
\label{eqn:WDVVa}
\end{eqnarray}
on $\overline{\cal{M}}_{0,a}$ with $a \geq 4$ (see, for example, Theorem 7.1 on page 599 in \cite{ACG}).
If $a=4$, then $x_5, \ldots, x_a$ should be dropped from the above equation.
Note that equations  \eqref{psi genus 0}, \eqref{psi genus 1}, and \eqref{eqn:WDVVa} hold for
any permutation of the set $\{x_1, \ldots, x_a \}$.

In the rest part of this section, we will give more details for the proofs of $B^2_{1, (2,1)}=0$ and $B^2_{1, (1,1,1)}=0$.
We will use Convention \ref{convention}
for notations of different types of legs. In particular, $W_i$ are the extra legs which will be forgotten
when applying $e_*$.

Equation \eqref{eqn:WDVVa} for the case $a=4$ is equivalent to say
$ <x_1 x_2 \gamma>_0 \, <\gamma^* x_3 x_4>_0$ is symmetric with respect to permutations of $\{x_1, \ldots, x_4\}$.
Repeatedly using this property, we can see that
\begin{equation} \label{eqn:g0p3Symm}
 <x_1 x_2 \gamma_1>_0 \, \bigg( \prod_{i=1}^k <\gamma_i^* x_{i+2} \gamma_{i+1}>_0 \bigg) <\gamma_{k+1}^* x_{k+3} x_{k+4}>_0
\end{equation}
is symmetric with respect to permutations of $\{x_1, \ldots, x_{k+4}\}$ for all $k \geq 0$.
We will use this property without mentioning it when combining like terms in all our calculations.

\subsection{Computing $B^2_{1, (2,1)}$}

Since there are only $2$ regular legs,  each $T \in \mathcal{A}_{1,2,2}$ can have at most $2$ top vertices.
Stability of $e'(T)$ implies that there can be only one top vertex.
By equation \eqref{eqn:pvrange}, non-trivial contributions to $B^2_{1, (2,1)}$ are given by
\begin{eqnarray*}
 B^2_{1, (2,1)}
&=& <V_1 V_2 \Psi^2 (U_1) \Psi (U_2)>_1
   -e_* \left( <V_1 V_2 \gamma>_0 <\gamma^* W_1 \Psi^2 (U_1) \Psi (U_2)>_1 \right) \\
&&  -e_* \left( <V_1 V_2 \Psi (U_2) \gamma>_0 <\gamma^* W_1 \Psi^2 (U_1)>_1 \right)  \\
&& -e_* \left( <V_1 V_2 \Psi^2 (\gamma)>_1 <\gamma^* W_1 W_2 W_3 \Psi^2 (U_1) \Psi (U_2)>_0 \right)\\
&& +e_* \left( <V_1 V_2 \gamma_1>_0 <\gamma^*_1 W_1 \Psi^2(\gamma_2)>_1 <\gamma^*_2 W_2 W_3 W_4 \Psi^2 (U_1) \Psi (U_2)>_0 \right) \\
&& +e_* \left( <V_1 V_2 \gamma_1>_0 <\gamma^*_1 W_1 \Psi (U_2) \gamma_2>_0 <\gamma^*_2 W_2 \Psi^2 (U_1)>_1 \right).
\end{eqnarray*}
Using equations \eqref{string equation} and \eqref{eqn:e*T} to get rid of $e_*$, we have
\begin{eqnarray}
 B^2_{1, (2,1)}
&=& <V_1 V_2 \Psi^2 (U_1) \Psi (U_2)>_1  - <V_1 V_2 \gamma>_0 <\gamma^* \Psi^2 (U_1) U_2>_1  \nonumber \\
&&  - <V_1 V_2 \gamma>_0 <\gamma^* \Psi (U_1) \Psi (U_2)>_1
        -3<V_1 V_2 \Psi^2 (\gamma)>_1 <\gamma^* U_1 U_2>_0  \nonumber \\
&&  +3<V_1 V_2 \gamma_1>_0 <\gamma^*_1 \Psi(\gamma_2)>_1 <\gamma^*_2 U_1 U_2>_0
  -<V_1 V_2 \Psi (U_2) \gamma>_0 <\gamma^* \Psi (U_1)>_1  \nonumber \\
&&    +<V_1 V_2 \gamma_1>_0 <\gamma^*_1 U_2 \gamma_2>_0 <\gamma^*_2 \Psi (U_1)>_1.
    \label{eqn:B1(21)}
\end{eqnarray}
Note that the last two terms in the right hand side of this equation cancel with each other
after applying equation \eqref{psi genus 0} with $a=4$ to $<V_1 V_2 \Psi (U_2) \gamma>_0$.
To compute the first 5 terms on the right hand side, we first derive some formulas for
higher powers of $\psi$-classes on moduli spaces of genus-0 and genus-1 curves.

Multiplying both sides of equation \eqref{psi genus 0} with $a=5$ by $\psi_2$ and
applying equation \eqref{psi genus 0} with $a=4$, we obtain
\begin{equation} \label{eqn:psi11M05}
<\Psi(x_1) \Psi(x_2) x_3 x_4 x_5>_0 = 2 <x_1 x_2 \gamma_1^*>_0 <\gamma_1 x_3 \gamma_2^*>_0 <\gamma_2 x_4 x_5>_0
\end{equation}
on $\overline{\mathcal{M}}_{0,5}$.
Multiplying both sides of equation \eqref{psi genus 1} with $a=4$ by $\psi_1 \psi_2$ and applying
equation \eqref{eqn:psi11M05} and  equation \eqref{psi genus 0} with $a=4$, we have
\begin{eqnarray}
<\Psi^2(x_1) \Psi(x_2) x_3 x_4>_1
&=& <x_1 x_3 \gamma_1^*>_0 \, <\gamma_1 x_4 \gamma_2^*>_0 \, <\gamma_2 \Psi(x_2)>_1 \nonumber \\
&& + 2 <x_1 x_2 \gamma_1^*>_0 <\gamma_1 x_3 \gamma_2^*>_0 <\gamma_2 x_4 \gamma_3^*>_0 \, <\gamma_3>_1 \nonumber \\
&& + \frac{1}{12} <\Psi(x_1) \Psi(x_2) x_3 x_4 \gamma^* \gamma>_0   \label{eqn:psi21M14-1}
\end{eqnarray}
on $\overline{\mathcal{M}}_{1,4}$.
Multiplying both sides of equation \eqref{psi genus 1} with $a=3$ by $\psi_1$ or $\psi_2$, we have
\begin{eqnarray}
<\Psi^2(x_1) x_2 x_3>_1
&=& <x_1 x_2 \gamma_1^*>_0 \, <\gamma_1 x_3 \gamma_2^*>_0 \, <\gamma_2 >_1 \nonumber \\
&& + \frac{1}{12} <\Psi(x_1) x_2 x_3 \gamma^* \gamma>_0   \label{eqn:psi2M13-1}
\end{eqnarray}
and
\begin{eqnarray}
<\Psi(x_1) \Psi(x_2) x_3>_1
&=& <x_1 x_3 \gamma^*>_0 \, <\gamma \Psi(x_2)>_1 \nonumber \\
&& + <x_1 x_2 \gamma_1^*> <\gamma_1 x_3 \gamma_2^*>_0 \, <\gamma_2 >_1 \nonumber \\
&& + \frac{1}{12} <x_1 \Psi(x_2) x_3 \gamma^* \gamma>_0   \label{eqn:psi11M13-1}
\end{eqnarray}
on $\overline{\mathcal{M}}_{1,3}$.

Applying equations \eqref{eqn:psi21M14-1}, \eqref{eqn:psi2M13-1}, \eqref{eqn:psi11M13-1}, and \eqref{psi genus 1} with $a=2$
to the first 5 terms on the right hand side of equation \eqref{eqn:B1(21)}, we can see immediately that all
terms containing genus-1 factors are cancelled and we have
\begin{eqnarray}
12 \, B^2_{1, (2,1)}
&=& <\Psi(U_1) \Psi(U_2) V_1 V_2 \gamma^* \gamma>_0 - <V_1 V_2 \gamma_1^*>_0 <\gamma_1 \Psi(U_1) U_2 \gamma_2^* \gamma_2>_0  \nonumber \\
&&  - <V_1 V_2 \gamma_1^*>_0 <\gamma_1 U_1 \Psi(U_2) \gamma_2^* \gamma_2>_0
      - 3 <U_1 U_2 \gamma_1^*>_0 <\Psi(\gamma_1) V_1 V_2 \gamma_2^* \gamma_2>_0
         \nonumber \\
&& + 3 <U_1 U_2 \gamma_1^*>_0 <\gamma_1 \gamma_2^* \gamma_2 \gamma_3^*>_0 <\gamma_3 V_1 V_2>_0.
    \label{eqn:B1(21)g0}
\end{eqnarray}

To compute the first term on the right hand side of the above equation, we multiply both sides of
equation \eqref{psi genus 0} with $a=6$ by $\psi_2$, then apply the formula with
\[x_1= U_1, x_2 = U_2, x_3=\gamma^*, x_4=\gamma, x_5=V_1, x_6=V_2\]
and use the symmetry for switching $\gamma^*$ and $\gamma$. We obtain the following formula
\begin{eqnarray}
&& <\Psi(U_1) \Psi(U_2) V_1 V_2 \gamma^* \gamma>_0 \nonumber \\
&=&  2 <U_1 \gamma^* \gamma_1^*>_0 <\gamma_1 \Psi(U_2) \gamma V_1 V_2>_0
   + 2 <U_1 U_2 \gamma_1^* >_0 <\gamma_1 \gamma^* \gamma_2^* >_0 <\gamma_2 \gamma V_1 V_2>_0  \nonumber \\
&& + <U_1 \gamma^* \gamma \gamma_1^* >_0 <\gamma_1 U_2 \gamma_2^* >_0 < \gamma_2 V_1 V_2 >_0 \nonumber \\
&&  + <U_1 \Psi(U_2) \gamma^* \gamma \gamma_1^* >_0 <\gamma_1 V_1 V_2>_0,
            \label{eqn:psi11M06-1}
\end{eqnarray}
where we have also used equation \eqref{psi genus 0} with $a=4$ to get rid of $\psi$-classes
in genus-$0$ factors with 4 insertions. Note that the last term on the right hand side of this equation
is cancelled with the third term on the right hand side of equation \eqref{eqn:B1(21)g0}.

After plugging \eqref{eqn:psi11M06-1} into equation \eqref{eqn:B1(21)g0}, applying
equation \eqref{psi genus 0} with $a=5$ to
\[ <\gamma_1 \Psi(U_2) \gamma V_1 V_2>_0 {\rm \,\,\, and \,\,\, }
 <\Psi(\gamma_1) V_1 V_2 \gamma_2^* \gamma_2>_0
 \]
by setting $x_4=V_1$ and $x_5 =V_2$,
applying equation \eqref{psi genus 0} with $a=5$ to
$<\gamma_1 \Psi(U_1) U_2 \gamma_2^* \gamma_2>_0$ by setting $x_4=\gamma_2^*$ and $x_5 =\gamma_2$,
we see immediately that coefficients of $<\gamma_2 \gamma V_1 V_2>_0$ add up to be $0$
and
\begin{equation} \label{eqn:BfV12}
 12 \, B^2_{1, (2,1)} = f <\gamma_2 V_1 V_2>_0,
 \end{equation}
where
\begin{eqnarray}
f
&=& 2 < U_1 \gamma^* \gamma_1^* >_0 <\gamma_1  U_2 \gamma \, \gamma_2^*  >_0
    + < U_2 \gamma_2^* \gamma_1^* >_0 <\gamma_1  U_1 \gamma^* \gamma >_0  \nonumber \\
&&  -  < U_1 \gamma_2^* \gamma_1^* >_0 <\gamma_1  U_2 \gamma^* \, \gamma  >_0
    - < U_1 U_2 \gamma_2^* \gamma_1^* >_0 <\gamma_1  \gamma^* \gamma >_0  \nonumber \\
&&   - < U_1 U_2 \gamma_1^* >_0 <\gamma_1  \gamma^* \gamma \, \gamma_2^*  >_0.
        \label{eqn:f}
\end{eqnarray}
Applying equation \eqref{eqn:WDVVa} with $a=5$ to the first two terms and last two terms separately, we
have
\begin{eqnarray}
&& < U_1 \gamma^* \gamma_1^* >_0 <\gamma_1  U_2 \gamma \, \gamma_2^*  >_0
    + < U_2 \gamma_2^* \gamma_1^* >_0 <\gamma_1  U_1 \gamma^* \gamma >_0 \nonumber \\
&=& < U_1 \gamma_2^* \gamma_1^* >_0 <\gamma_1  U_2 \gamma \, \gamma^*  >_0
    + < U_2 \gamma^* \gamma_1^* >_0 <\gamma_1  U_1 \gamma_2^* \gamma >_0
    \label{eqn:WDVV5-1}
\end{eqnarray}
and
\begin{eqnarray}
&& < U_1 U_2 \gamma_2^* \gamma_1^* >_0 <\gamma_1  \gamma^* \gamma >_0
    + < U_1 U_2 \gamma_1^* >_0 <\gamma_1  \gamma^* \gamma \, \gamma_2^*  >_0 \nonumber \\
&=& < U_1 \gamma^*  \gamma_2^* \gamma_1^* >_0 <\gamma_1  U_2 \gamma >_0
    + < U_1 \gamma^* \gamma_1^* >_0 <\gamma_1   U_2 \gamma \, \gamma_2^*  >_0.
    \label{eqn:WDVV5-2}
\end{eqnarray}

Plugging formulas \eqref{eqn:WDVV5-1} and \eqref{eqn:WDVV5-2} into equation \eqref{eqn:f},
we see immediately that $f=0$. Hence by equation \eqref{eqn:BfV12}, we have
proved the following
\begin{prop}\label{B21=0}
$B^2_{1,(2,1)}=0$.
\end{prop}

\subsection{Computing $B^2_{1, (1,1,1)}$}

Let ${\mathcal C}_3$ be the set of cyclic permutations of $(1, 2, 3)$, i.e.
\[ {\mathcal C}_3 := \{(1, 2, 3), (2, 3, 1), (3, 1, 2) \}. \]
After getting rid of $e_*$ and combining like terms, $B^2_{1,(1,1,1)}$ can be written in the following form
\begin{equation}
 B^2_{1,(1,1,1)} = H+<V_1 V_2 \gamma>_0 I, \label{eqn:BHI}
\end{equation}
where
\begin{eqnarray}
H &:=& <V_1 V_2 \Psi(U_1) \Psi(U_2) \Psi(U_3)>_1 -6 <V_1 V_2 \Psi^2 (\gamma)>_1 <\gamma^* U_1 U_2 U_3>_0
   \nonumber \\
&& -2 \sum_{(i,j,k) \in {\mathcal C}_3}
                 <V_1 V_2 \Psi (U_i) \Psi (\gamma)>_1 <\gamma^* U_j U_k>_0
    \nonumber \\
&&
 + 6 <V_1 V_2 \Psi(\gamma_1)>_1 <\gamma^*_1 U_1 \gamma_2>_0 <\gamma^*_2 U_2 U_3>_0, \label{eqn:H}
\end{eqnarray}
 and
\begin{eqnarray}
I &:=& 6<\gamma^* \Psi (\gamma_1)>_1 <\gamma^*_1 U_1 U_2 U_3>_0
-6<\gamma^* \gamma_1>_1 <\gamma^*_1 U_1 \gamma_2>_0 <\gamma^*_2 U_2 U_3>_0
        \nonumber \\
&&
+\sum_{(i,j,k) \in {\mathcal C}_3}
\big(
 -<\gamma^* U_i \Psi (U_j) \Psi (U_k)>_1  \nonumber \\
&& \hspace{20pt} + 2<\gamma^* \Psi (U_i) \gamma_1>_1 <\gamma^*_1 U_j U_k>_0
+2<\gamma^* U_i \Psi (\gamma_1)>_1 <\gamma^*_1 U_j U_k>_0
   \big). \label{eqn:I}
\end{eqnarray}

As in the proof of Proposition \ref{B21=0}, we first use equations \eqref{psi genus 0} and \eqref{psi genus 1} to get rid of $\psi$-classes in $B^2_{1,(1,1,1)}$, then use equation \eqref{eqn:WDVVa} to prove that it is equal to 0. However, the computation for $B^2_{1,(1,1,1)}$ is much more complicated than that for $B^2_{1,(2,1)}$. To simplify computations, we will adopt the following strategy: We do not remove all $\psi$-classes for each term in $B^2_{1,(1,1,1)}$ at the beginning since this will result in a very complicated expression. Instead, we consider terms with higher degrees of $\psi$-classes on a genus-1 vertex and remove one $\psi$-class first, then combine them with other terms to simplify the expression before removing more $\psi$-classes.

Note that there is no ambiguity in applying equation \eqref{psi genus 1} when it is clear which $\psi$-class should be removed on a genus-1 vertex. But when applying equation \eqref{psi genus 0} to remove a $\psi$-class on a genus-0 vertex, different choices of $x_{a-1}$ and $x_a$ will produce expressions which look quite differently if $a>4$. This could cause
a lot of troubles when combining like terms. To reduce such troubles, we will obey the following rules when applying equation \eqref{psi genus 0} to a genus-0 vertex with at least one $\psi$-class: If $V_1$ and $V_2$ are both attached to this vertex, we always choose them to be  $x_{a-1}$ and $x_a$.
If it is possible, we should choose $x_{a-1}$ and $x_a$ to be two legs, i.e. $V_i$ or $U_j$ for some $i$ and $j$.
If both half edges of an edge are attached to this vertex, we should avoid choose them to be either $x_{a-1}$ or $x_a$ if possible.

Recall that $H$ and $I$ are defined by equations \eqref{eqn:H} and \eqref{eqn:I} respectively. We will compute $H$ and $I$ separately. We will show that after getting rid of all $\psi$-classes, the genus-1 parts of $H$ and $I$ are both equal to $0$, and the sum of purely genus-0 contributions of $H$ and $I$ to $B^2_{1,(1,1,1)}$ is also equal to $0$.

\subsubsection{Computing $H$}

Consider the term $<V_1 V_2 \Psi(U_1) \Psi(U_2) \Psi(U_3)>_1$ in $H$. We first remove
the $\psi$-class associated with $U_1$ in this term using
equation \eqref{psi genus 1}.
In the resulting expression, if there are genus-1 vertices containing both $\psi(U_2)$ and $\psi(U_3)$, we
remove the $\psi$-class associated with $U_2$ using
equation \eqref{psi genus 1}.
We then remove all $\psi$-classes associated to half edges $\gamma$ or $\gamma_1$ in other terms in $H$.
Putting all resulting expressions together and combing like terms first, we then use
equations \eqref{psi genus 0} and \eqref{psi genus 1} to remove all $\psi$-classes
in all terms which have a genus-1 factor. We then combine all terms which can factor out a common genus-1 factor.
Using symmetry for expression \eqref{eqn:g0p3Symm}, we can show that coefficients of genus-1 vertices with at least two half edges are all equal to $0$, and $H$ has the following form
\begin{equation} \label{eqn:DefH}
H = 2 <\gamma_1>_1 H_1 + H_0,
\end{equation}
where
\begin{eqnarray}
 H_1
&:=& 3<\gamma^*_1 U_2 \gamma_2>_0 <\gamma^*_2 U_3 \gamma_3>_0 <\gamma^*_3 U_1 V_1 V_2>_0 \nonumber \\
&&    + \big\{ <\gamma^*_1 V_1 V_2 \gamma_2>_0 <\gamma^*_2 U_1 \gamma_3>_0
      -2<\gamma^*_1 U_1 \gamma_2>_0 <\gamma^*_2 V_1 V_2 \gamma_3>_0  \big\} <\gamma^*_3 U_2 U_3>_0 \nonumber \\
&&- \sum_{i=2}^3
    <\gamma^*_1 U_i \gamma_2>_0 <\gamma^*_2 V_1 V_2 \gamma_3>_0 <\gamma^*_3 U_1 U_{5-i}>_0
   \nonumber \\
&& + \sum_{i=1}^2 \big( -<\gamma^*_1 U_1 V_i \gamma_2>_0 <\gamma^*_2 U_2 \gamma_3>_0 <\gamma^*_3 U_3 V_{3-i}>_0
           \nonumber \\
&& \hspace{60pt} +<\gamma^*_1 U_3 \gamma_2>_0 <\gamma^*_2 U_2 V_i \gamma_3>_0 <\gamma^*_3 U_1 V_{3-i}>_0
      \nonumber \\
&&\hspace{60pt} +<\gamma^*_1 U_3 \gamma_2 \gamma_3>_0 <\gamma^*_2 U_1 V_i>_0 <\gamma^*_3 U_2 V_{3-i}>_0
   \big) \nonumber \\
&&+<V_1 V_2 \gamma>_0
\big( <\gamma^* U_3 \gamma_2>_0 <\gamma^*_2 U_1 U_2 \gamma_1^*>_0
                -3<\gamma^* \gamma_1^* \gamma_2>_0 <\gamma^*_2 U_1 U_2 U_3>_0 \nonumber \\
&&\hspace{90pt}  - \sum_{i=1}^2 <\gamma^* \gamma_1^* U_i \gamma_2>_0 <\gamma^*_2 U_{3-i} U_3>_0
        \nonumber \\
&& \hspace{90pt} + \sum_{i=2}^3 <\gamma^* U_1 U_i \gamma_2>_0 <\gamma^*_2 U_{5-i} \gamma_1^*>_0
\big),
\label{second formula throw genus 0}
\end{eqnarray}
and
\begin{eqnarray}
12 H_0
&:=& <U_1 \Psi (U_2) \Psi (U_3) V_1 V_2 \gamma \gamma^*>_0
 \nonumber \\
&&-6<U_1 U_2 U_3 \gamma_1>_0 \big( <\gamma^*_1 V_1 V_2 \gamma_2>_0 <\gamma^*_2 \gamma \gamma^*>_0
   + \sum_{i=1}^2 <\gamma^*_1 V_i \gamma_2>_0 <\gamma^*_2 V_{3-i} \gamma \gamma^*>_0 \big)
    \nonumber \\
&&  -2 \sum_{(i,j,k) \in {\mathcal C}_3}
          <U_j U_k \gamma_1>_0 <\gamma^*_1 \Psi (U_i) V_1 V_2 \gamma \gamma^*>_0  \nonumber \\
&& +<U_1 V_1 V_2 \gamma_1>_0 \big( <\gamma^*_1 U_2 \Psi (U_3) \gamma \gamma^*>_0
                                  + 2 <\gamma^*_1 U_2 \gamma_2>_0 <\gamma^*_2 U_3 \gamma \gamma^*>_0 \nonumber \\
&& \hspace{120pt}                  + <\gamma^*_1 U_3 \gamma_2>_0 <\gamma^*_2 U_2 \gamma \gamma^*>_0  \big)
             \nonumber \\
&& + <V_1 V_2 \gamma_1>_0 \sum_{i=2}^3 <\gamma^*_1 U_1 U_{i} \gamma_2>_0 <\gamma^*_2 U_{5-i} \gamma \gamma^*>_0
                                \nonumber \\
&& + \sum_{i=1}^2 <U_1 V_i \gamma_1>_0 \big( <\gamma^*_1 U_2 \Psi (U_3) V_{3-i} \gamma \gamma^*>_0
                         + <\gamma^*_1 U_3 \gamma \gamma^* \gamma_2>_0 <\gamma^*_2 U_2 V_{3-i}>_0  \nonumber \\
&& \hspace{120pt}
                             + <\gamma^*_1 U_2 V_{3-i} \gamma_2>_0 <\gamma^*_2 U_3 \gamma \gamma^*>_0 \big)
                             \nonumber \\
&& - 2 <U_2 U_3 \gamma_1>_0 \big(<\gamma^*_1 V_1 V_2 \gamma_2>_0 <\gamma^*_2 U_1 \gamma \gamma^*>_0
            + \sum_{i=1}^2 <\gamma^*_1 V_i \gamma_2>_0 <\gamma^*_2 U_1 V_{3-i} \gamma \gamma^*>_0 \big)
                                            \nonumber \\
&& - <U_1 U_3 \gamma_1>_0 \big(<\gamma^*_1 V_1 V_2 \gamma_2>_0 <\gamma^*_2 U_2 \gamma \gamma^*>_0
      + \sum_{i=1}^2 <\gamma^*_1 V_i \gamma_2>_0 <\gamma^*_2 U_2 V_{3-i} \gamma \gamma^*>_0 \big) \nonumber \\
&&- <U_1 U_2 \gamma_1>_0 <\gamma^*_1 V_1 V_2 \gamma_2>_0 <\gamma^*_2 U_3 \gamma \gamma^*>_0 \nonumber \\
&&+ 6 <V_1 V_2 \gamma \gamma^* \gamma_1>_0 <\gamma^*_1 U_1 \gamma_2>_0 <\gamma^*_2 U_2 U_3>_0.
                \label{second thrown formula}
\end{eqnarray}

We will repeatedly apply equation \eqref{eqn:WDVVa} to show $H_1=0$. For the purpose of presentation,
it is convenient to introduce the following formal notation:
\begin{equation} \label{eqn:partial}
\partial_x \big( \prod_{i=1}^{k}<y_{a_{i-1}} \cdots y_{a_i}>_{g_i}\big) := \sum_{j=1}^{k} <x y_{a_{j-1}} \cdots y_{a_j}>_{g_j} \prod_{1 \leq i \leq k \atop i \neq j}<y_{a_{i-1}} \cdots y_{a_i}>_{g_i},
\end{equation}
 where $x$ and $y_{a_i}$ are half edges. Using this notation, equation \eqref{eqn:WDVVa} can be written in a
 concise form
\begin{equation} \label{eqn:WDVVaP}
\partial_{x_5} \ldots \partial_{x_a} \big( <x_1 x_2 \gamma>_0 <\gamma^* x_3 x_4>_0
   -  <x_1 x_3 \gamma>_0 <\gamma^* x_2 x_4>_0 \big) =0.
\end{equation}

We are now ready to prove
\begin{lem}\label{H1=0}
$H_1=0$.
\end{lem}

\noindent
{\bf Proof}:
It is straightforward to check the following equality by first expanding the right hand side using
definition of "$\partial$" given by equation \eqref{eqn:partial} and then combining like terms:
\begin{eqnarray}
H_1
&=& \partial_{V_1} \big( <\gamma^*_1 V_2 \gamma_2>_0 <\gamma^*_2 U_3 \gamma_3>_0
                        -<\gamma^*_1 U_3 \gamma_2>_0 <\gamma^*_2 V_2 \gamma_3>_0 \big)
                        <\gamma^*_3 U_1 U_2>_0 \nonumber \\
&&+ \sum_{i=1}^2 \partial_{V_i} \big( <\gamma^*_1\gamma_2\gamma_3>_0 <\gamma^*_2 U_1 V_{3-i}>_0
            -<\gamma^*_1 U_1 \gamma_2>_0 <\gamma^*_2 V_{3-i} \gamma_3>_0 \big)
            <\gamma^*_3 U_2 U_3>_0 \nonumber \\
&& +\partial_{\gamma_3} \big( <\gamma^*_1 U_3 \gamma_2>_0 <\gamma^*_2 U_2 V_2>_0
             -<\gamma^*_1 V_2 \gamma_2>_0 <\gamma^*_2 U_2 U_3>_0 \big)
             <\gamma^*_3 U_1 V_1>_0 \nonumber \\
&&+ \big\{ \partial_{\gamma_1^*} \big( <U_3 \gamma_2 \gamma_3>_0 <\gamma^*_2 U_2 V_1>_0
       -<V_1 \gamma_2 \gamma_3>_0 <\gamma^*_2 U_2 U_3>_0 \big)
       \nonumber \\
&&\hspace{20pt} +\partial_{V_1} \big( <\gamma^*_1 U_3 \gamma_2>_0 <\gamma^*_2 U_2 \gamma_3>_0
            -<\gamma^*_1 U_2 \gamma_2>_0 <\gamma^*_2 U_3 \gamma_3>_0 \big)
            \big\} <\gamma^*_3 U_1 V_2>_0 \nonumber \\
&&+<\gamma^*_1 U_2 \gamma_2>_0
    \partial_{V_1} \big( <\gamma^*_2 U_3 \gamma_3>_0 <\gamma^*_3 U_1 V_2>_0
                        -<\gamma^*_2 V_2 \gamma_3>_0 <\gamma^*_3 U_1 U_3>_0 \big) \nonumber \\
&&+<\gamma^*_1 V_2 \gamma_2>_0
    \partial_{U_3} \big( <\gamma^*_2 U_2 \gamma_3>_0 <\gamma^*_3 U_1 V_1>_0
                         -<\gamma^*_2 V_1 \gamma_3>_0 <\gamma^*_3 U_1 U_2>_0 \big) \nonumber \\
&&+<V_1 V_2 \gamma>_0 \big\{ \partial_{U_1} \big( <\gamma^*_1 U_2 \gamma_2>_0 <\gamma^*_2 U_3 \gamma^*>_0
                     -<\gamma^* \gamma^*_1 \gamma_2>_0 <\gamma^*_2 U_2 U_3>_0 \big) \nonumber \\
&& \hspace{80pt} +\partial_{U_2} \big( <\gamma^* U_1 \gamma_2>_0 <\gamma^*_2 U_3 \gamma^*_1>_0
                     -<\gamma^* \gamma^*_1 \gamma_2>_0 <\gamma^*_2 U_1 U_3>_0 \big) \big\}. \nonumber
\end{eqnarray}
The right hand side of this equality is equal to 0 due to equation \eqref{eqn:WDVVaP}.
The lemma is thus proved.
$\Box$

\subsubsection{Computing $I$}

Recall that $I$ is defined by equation \eqref{eqn:I}.
We first remove one $\psi$-class on every genus-1 vertex with two $\psi$-classes using equation \eqref{psi genus 1}.
To keep the symmetry of $U_i$, we remove the $\psi$-class associated to $U_1$
for $<\gamma^* \Psi (U_1) \Psi (U_2) U_3>_1$,
remove the $\psi$-class associated to $U_2$ for $<\gamma^* U_1 \Psi (U_2) \Psi (U_3)>_1$,
and remove the $\psi$-class associated to $U_3$ for $<\gamma^* \Psi (U_1) U_2 \Psi (U_3)>_1$.
After combining like terms, we remove $\psi$-classes from all terms which has a genus-1 factor
using equations \eqref{psi genus 0} and \eqref{psi genus 1}. We then obtain
\begin{equation} \label{eqn:I=1+0}
I \, = \, <\gamma_1>_1 I_1 + I_0,
\end{equation}
where
\begin{eqnarray}
I_1 &=& \sum_{(i,j,k) \in {\mathcal C}_3}
       \big( <\gamma^* \gamma^*_1 U_i \gamma_2>_0 <\gamma^*_2 U_j U_k>_0
              - <\gamma^* U_i \gamma_2>_0 <\gamma^*_2 U_j U_k \gamma^*_1>_0
              \nonumber \\
&& \hspace{60pt} - <\gamma^*_1 U_i \gamma_2>_0 <\gamma^*_2 U_j U_k \gamma^*>_0
               \big) \nonumber \\
&& +3<\gamma^* \gamma^*_1 \gamma_2>_0 <\gamma^*_2 U_1 U_2 U_3>_0, \label{fourth formula throw genus 0}
\end{eqnarray}
and
\begin{eqnarray}
12 I_0 &=& \sum_{(i,j,k) \in {\mathcal C}_3}
       \big( - <\gamma^*   \Psi (U_i) U_j U_k \gamma_1 \gamma^*_1>_0
        + 3 <\gamma^* U_i \gamma_1 \gamma_2 \gamma_2^*>_0 <\gamma^*_1 U_j U_k>_0
     \nonumber \\
&&   \hspace{30pt} - <\gamma^* U_i \gamma_1>_0 <\gamma_1^* U_j U_k \gamma_2 \gamma_2^*>_0
     - <\gamma^* U_j U_k \gamma_1>_0 <\gamma_1^* U_i \gamma_2 \gamma_2^*>_0 \big)
      \nonumber \\
&&  + 6<\gamma^* \gamma_1 \gamma_2 \gamma_2^*>_0 <\gamma^*_1 U_1 U_2 U_3>_0.
\label{third thrown formula}
\end{eqnarray}

We have the following
\begin{lem}\label{I1=0}
$I_1=0$.
\end{lem}

\noindent
{\bf Proof}:
It is straightforward to check the following equality
by expanding the right hand side and combining like terms:
\begin{eqnarray}
I_1 &=& \sum_{(i,j,k) \in {\mathcal C}_3}
        \partial_{U_i} \big( <\gamma^* \gamma^*_1 \gamma_2>_0 <\gamma^*_2 U_j U_k>_0
                 -<\gamma^* U_j \gamma_2>_0 <\gamma^*_2 U_k \gamma^*_1>_0 \big) . \nonumber
\end{eqnarray}
The right hand side of this equation is $0$ by equation \eqref{eqn:WDVVaP}.
The lemma is thus proved.$\Box$

\subsubsection{Combining $H$ and $I$}

We are now ready to prove the following
\begin{prop} \label{B111=0}
$B^2_{1,(1,1,1)}=0$.
\end{prop}

\noindent
{\bf Proof}:
By equations \eqref{eqn:BHI}, \eqref{eqn:DefH}, \eqref{eqn:I=1+0}, and Lemmas \ref{H1=0} and \ref{I1=0}, we have
\begin{equation} \label{eqn:B111g0}
B^2_{1,(1,1,1)} = H_0 +<V_1 V_2 \gamma>_0 I_0,
\end{equation}
where $H_0$ and $I_0$ are defined by equations \eqref{second thrown formula} and \eqref{third thrown formula} respectively.
Hence we only need to show that the right hand side of equation \eqref{eqn:B111g0} is equal to 0.

After eliminating $\psi$-classes in $H_0$ and $I_0$ using equation \eqref{psi genus 0}  and combining like terms, we obtain
\begin{eqnarray}
&& 6 \big( H_0 +<V_1 V_2 \gamma>_0 I_0 \big) \nonumber \\
&=&3 \big( <V_1 V_2 \gamma_1>_0 <\gamma^*_1  \gamma_2 \gamma^*_2 \gamma_3>_0
           - <\gamma_1 \gamma^*_1 \gamma_2>_0 <\gamma_2^* V_1 V_2 \gamma_3>_0
           \nonumber \\
&& \hspace{20pt}  - \sum_{i=1}^2 <\gamma_1 \gamma^*_1 V_i \gamma_2>_0 <\gamma_2^* V_{3-i} \gamma_3>_0
           \big) <\gamma^*_3 U_1 U_2 U_3>_0
        \nonumber \\
&& + \big( <\gamma_1 \gamma^*_1 V_1 V_2 \gamma_2>_0 <\gamma_2^* U_3 \gamma_3>_0
            -<U_3 \gamma_1 \gamma_2>_0 <\gamma^*_1 \gamma_2^* V_1 V_2 \gamma_3>_0
            \big) <\gamma^*_3 U_1 U_2>_0
       \nonumber \\
&&- \big( <U_2 \gamma_1 \gamma_2>_0 <\gamma^*_1 \gamma_2^* V_1 V_2 \gamma_3>_0
      + 2<V_1 V_2 \gamma_1 \gamma_2>_0 <\gamma^*_1 \gamma^*_2 U_2 \gamma_3> \big) <\gamma^*_3 U_1 U_3>_0 \nonumber \\
&&+ \big( 3<U_1 V_1 V_2 \gamma_1 \gamma_2>_0 <\gamma^*_1 \gamma_2^* \gamma_3>_0
           -2<U_1 \gamma_1 \gamma_2>_0 <\gamma^*_1 \gamma_2^* V_1 V_2 \gamma_3>_0  \nonumber \\
&& \hspace{10pt}
          -2<V_1 V_2 \gamma_1 \gamma_2>_0 <\gamma^*_1 \gamma^*_2 U_1 \gamma_3>
          - \sum_{i=1}^2 <\gamma_1 \gamma^*_1 U_1 V_i \gamma_2>_0 <\gamma_2^* V_{3-i} \gamma_3>_0
          \big) <\gamma^*_3 U_2 U_3>_0
        \nonumber \\
&&+ \big(2 <\gamma_1 \gamma^*_1 U_3 \gamma_2>_0 <\gamma_2^* U_2 \gamma_3>_0
          + <\gamma_1 \gamma^*_1 U_2 \gamma_2>_0 <\gamma_2^* U_3 \gamma_3>_0
           \nonumber \\
&& \hspace{10pt} +<U_2 \gamma_1 \gamma_2>_0 <\gamma^*_1 \gamma^*_2 U_3 \gamma_3>_0
           +2<U_3 \gamma_1 \gamma_2>_0 <\gamma^*_1 \gamma^*_2 U_2 \gamma_3>_0
           \big) <\gamma_3^* U_1 V_1 V_2>_0 \nonumber \\
&&+ \sum_{i=1}^2 \big( <\gamma_3 U_3 \gamma_2 \gamma^*_2 \gamma_1>_0 <\gamma^*_1 U_2 V_{3-i}>_0
         +  <\gamma_1 \gamma^*_1 U_3 \gamma_2>_0 <\gamma_2^* U_2 V_{3-i} \gamma_3>_0
     \nonumber \\
&&  \hspace{20pt} + <U_3 \gamma_1 \gamma_2>_0 <\gamma^*_1 \gamma_2^* U_2 V_{3-i} \gamma_3>_0
                          +<U_2 V_{3-i} \gamma_1 \gamma_2>_0 <\gamma^*_1 \gamma^*_2 U_3 \gamma_3>_0
                          \big) <\gamma^*_3 U_1 V_i>_0  \nonumber \\
&& + \big(2<U_2 \gamma_1 \gamma_2>_0 <\gamma^*_1 U_1 U_3 \gamma_3>_0
  +4<U_3 \gamma_1 \gamma_2>_0 <\gamma^*_1 U_1 U_2 \gamma_3>_0 \nonumber \\
&& \hspace{10pt} - 2<\gamma_1 \gamma^*_1 U_1 \gamma_2>_0  <\gamma_3 U_2 U_3>_0
        - \sum_{i=2}^3 <\gamma_1 \gamma^*_1 U_i \gamma_2>_0 <\gamma_3 U_1 U_{5-i} >_0
            \big) <\gamma_2^* \gamma^*_3 V_1 V_2 >_0  \nonumber \\
&& + \big(  <U_1 U_2 \gamma_1>_0 <\gamma^*_1 U_3 \gamma_2 \gamma^*_2 \gamma_3>_0
            - \sum_{i=1}^2 <\gamma_1 \gamma^*_1 U_i U_3 \gamma_2>_0 <\gamma_2^* U_{3-i} \gamma_3>_0 \nonumber \\
&& \hspace{10pt}  - <\gamma_1 \gamma^*_1 U_1 \gamma_2>_0 <\gamma_2^* U_2 U_3 \gamma_3>_0
                -<U_1 \gamma_1 \gamma_2>_0 <\gamma^*_1 \gamma_2^* U_2 U_3 \gamma_3>_0 \nonumber \\
&&\hspace{10pt} + <U_1 U_2 \gamma_1 \gamma_2>_0 <\gamma^*_1 \gamma^*_2 U_3 \gamma_3>_0
                - \sum_{i=1}^2 <U_i U_3 \gamma_1 \gamma_2>_0 <\gamma^*_1 \gamma^*_2 U_{3-i} \gamma_3>
            \big) <\gamma^*_3 V_1 V_2>_0.
       \nonumber
\end{eqnarray}

It is straightforward to check the following equality
by expanding the right hand side and combining like terms:
\begin{eqnarray}
&& 6 \big( H_0 +<V_1 V_2 \gamma>_0 I_0 \big) \nonumber \\
&=& \big\{ 3\partial_{V_1} \big( <V_2 \gamma_1 \gamma_2>_0 <\gamma^*_1 \gamma^*_2 \gamma_3>_0
                                 -<\gamma_1 \gamma^*_1 \gamma_2>_0 <\gamma^*_2 V_2 \gamma_3>_0 \big) \nonumber \\
&& \hspace{10pt} +3\partial_{\gamma^*_1} \big(<\gamma_1 \gamma_2 \gamma_3>_0 <\gamma^*_2 V_1 V_2>_0
                            -<V_2 \gamma_1 \gamma_2>_0 <\gamma^*_2 V_1 \gamma_3>_0 \big) \nonumber \\
&&\hspace{10pt} + \partial_{\gamma_3} \big( <V_1 \gamma_1 \gamma_2>_0 <\gamma^*_1 \gamma^*_2 V_2>_0
                          -<V_1 V_2 \gamma_1>_0 <\gamma^*_1 \gamma_2 \gamma^*_2>_0 \big)
        \big\} <\gamma^*_3 U_1 U_2 U_3>_0 \nonumber \\
&&+\sum_{i=2}^3 \partial_{V_1} \big( <U_1 V_2 \gamma_1>_0 <\gamma^*_1 \gamma^*_3 U_{5-i}>_0
       -<U_1 U_{5-i} \gamma_1>_0 <\gamma^*_1 \gamma^*_3 V_2>_0 \big)
       <U_i \gamma_2 \gamma^*_2 \gamma_3>_0 \nonumber \\
&&+ \sum_{i=2}^3 \partial_{U_i} \big( <\gamma_1 \gamma^*_1 \gamma_2>_0 <\gamma^*_2 V_2 \gamma_3>_0
       -<V_2 \gamma_1 \gamma_2>_0 <\gamma^*_1 \gamma^*_2 \gamma_3>_0 \big) <\gamma^*_3 U_1 U_{5-i} V_1>_0 \nonumber \\
&&+ \sum_{i=1}^2 \partial_{\gamma_2} \big(<U_1 \gamma_1 \gamma_3>_0 <\gamma^*_1 \gamma^*_3 V_i>_0
             -<U_1 V_i \gamma_1>_0 <\gamma^*_1 \gamma_3 \gamma^*_3>_0 \big)
                          <\gamma^*_2 U_2 U_3 V_{3-i}>_0 \nonumber \\
&&+\partial_{U_3} \big( <\gamma_1 \gamma^*_1 \gamma_2>_0 <\gamma^*_2 U_2 \gamma_3>_0
       -<U_2 \gamma_1 \gamma_2>_0 <\gamma^*_1 \gamma^*_2 \gamma_3>_0 \big)
       <\gamma^*_3 U_1 V_1 V_2>_0 \nonumber \\
&& + \big\{ \partial_{V_2}\partial_{\gamma_1^*} \big( <V_1 \gamma_1 \gamma_2>_0 <\gamma^*_2 U_1 \gamma_3>_0
                                                 -<U_1 \gamma_1 \gamma_2>_0 <\gamma^*_2 V_1 \gamma_3>_0 \big) \nonumber \\
&& \hspace{20pt} +\partial_{U_1}\partial_{V_1} \big( <V_2 \gamma_1 \gamma_2>_0 <\gamma^*_1 \gamma^*_2 \gamma_3>_0
                             -<\gamma_1 \gamma^*_1 \gamma_2>_0 <\gamma^*_2 V_2 \gamma_3>_0 \big)
      \big\} <\gamma^*_3 U_2 U_3>_0 \nonumber \\
&&+\partial_{U_3}\partial_{\gamma_3} \big( <\gamma_1 \gamma^*_1 \gamma_2>_0 <\gamma^*_2 U_2 V_2>_0
     -<U_2 \gamma_1 \gamma_2>_0 <\gamma^*_1 \gamma^*_2 V_2>_0 \big) <\gamma^*_3 U_1 V_1>_0 \nonumber \\
&&+\partial_{U_3}\partial_{\gamma_1} \big( <U_2 V_1 \gamma_2>_0 <\gamma^*_1 \gamma^*_2 \gamma_3>_0
                -<U_2 \gamma_1^* \gamma_2>_0 <\gamma^*_2 V_1 \gamma_3>_0 \big) <\gamma^*_3 U_1 V_2>_0 \nonumber \\
&&+\partial_{U_3}\partial_{\gamma_1} \big( <U_1 U_2 \gamma_2>_0 <\gamma^*_1 \gamma^*_2 \gamma_3>_0
       -<U_1 \gamma_1^* \gamma_2>_0 <\gamma^*_2 U_2 \gamma_3>_0 \big) <\gamma^*_3 V_1 V_2>_0 \nonumber \\
&&+<U_2 U_3 \gamma_1 \gamma_2>_0 \big\{ \partial_{\gamma^*_1} \big( <\gamma^*_2 V_2 \gamma_3>_0 <\gamma^*_3 U_1 V_1>_0
           -<\gamma^*_2 U_1 \gamma_3>_0 <\gamma_3^* V_1 V_2>_0 \big) \nonumber \\
&& \hspace{100pt} + \partial_{V_1} \big( <\gamma^*_1 \gamma^*_2 \gamma_3>_0 <\gamma_3^* U_1 V_2>_0
             -<\gamma^*_1 U_1 \gamma_3>_0 <\gamma^*_2 \gamma_3^* V_2>_0 \big) \big\}\nonumber \\
&&+<U_1 \gamma_1 \gamma_2>_0 \partial_{V_1} \partial_{V_2} \big( <\gamma^*_1 U_2 \gamma_3>_0 <\gamma^*_2 \gamma^*_3 U_3>_0
                        -<\gamma^*_1 \gamma^*_2 \gamma_3>_0 <\gamma^*_3 U_2 U_3>_0 \big) \nonumber \\
&&+<U_2 \gamma_1 \gamma_2>_0 \big\{ \partial_{V_1}\partial_{\gamma_1^*}
                  \big( <\gamma^*_2 U_3 \gamma_3>_0 <\gamma^*_3 U_1 V_2>_0
                        -<\gamma^*_2 V_2 \gamma_3>_0 <\gamma^*_3 U_1 U_3>_0 \big) \nonumber \\
&&\hspace{90pt} + \partial_{U_3} \partial_{V_2} \big( <\gamma^*_1 \gamma^*_2 \gamma_3>_0 <\gamma^*_3 U_1 V_1>_0
                          -<\gamma^*_1 U_1 \gamma_3>_0 <\gamma^*_2 \gamma^*_3 V_1>_0 \big) \big\} \nonumber \\
&&+<U_3 \gamma_1 \gamma_2>_0
      \big\{ \partial_{V_1} \partial_{\gamma^*_1} \big( <\gamma^*_2 U_2 \gamma_3>_0 <\gamma^*_3 U_1 V_2>_0
                                        -<\gamma^*_2 V_2 \gamma_3>_0 <\gamma^*_3 U_1 U_2>_0 \big) \nonumber \\
&&\hspace{90pt} + \partial_{U_2} \partial_{V_2} \big( <\gamma^*_1 \gamma^*_2 \gamma_3>_0 <\gamma^*_3 U_1 V_1>_0
                          -<\gamma^*_1 U_1 \gamma_3>_0 <\gamma^*_2 \gamma^*_3 V_1>_0 \big)
             \big\}\nonumber \\
&&+ \sum_{i=2}^3 <U_i V_2 \gamma_1 \gamma_2>_0 \partial_{U_{5-i}}
        \big( <\gamma^*_1 \gamma^*_2 \gamma_3>_0 <\gamma_3^* U_1 V_1>_0
           -<\gamma^*_1 U_1 \gamma_3>_0 <\gamma^*_2 \gamma_3^* V_1>_0 \big) \nonumber \\
&&+<V_1 V_2 \gamma_1 \gamma_2>_0 \big\{
    4 \partial_{U_1} \big( <\gamma^*_1 U_2 \gamma_3>_0 <\gamma^*_2 \gamma_3^* U_3>_0
                    -<\gamma^*_1 \gamma^*_2 \gamma_3>_0 <\gamma_3^* U_2 U_3>_0 \big) \nonumber \\
&& \hspace{100pt}  +2 \partial_{U_2} \big( <\gamma^*_1 U_1 \gamma_3>_0 <\gamma^*_2 \gamma_3^* U_3>_0
              -<\gamma^*_1 \gamma^*_2 \gamma_3>_0 <\gamma_3^* U_1 U_3>_0 \big)
            \big\} \nonumber \\
&&+<V_1 \gamma_1 \gamma_2>_0
       \partial_{U_1} \partial_{V_2} \big( <\gamma^*_1 U_2 \gamma_3>_0 <\gamma^*_2 \gamma^*_3 U_3>_0
                                  -<\gamma^*_1 \gamma^*_2 \gamma_3>_0 <\gamma^*_3 U_2 U_3>_0 \big) \nonumber \\
&&+<V_2 \gamma_1 \gamma_2>_0 \big\{
            \partial_{U_2} \partial_{U_3} \big( <\gamma^*_1 \gamma^*_2 \gamma_3>_0 <\gamma^*_3 U_1 V_1>_0
                          -<\gamma^*_1 U_1 \gamma_3>_0 <\gamma^*_2 \gamma^*_3 V_1>_0 \big) \nonumber \\
&& \hspace{90pt}   + \partial_{U_1} \partial_{V_1} \big( <\gamma^*_1 U_2 \gamma_3>_0 <\gamma^*_2 \gamma^*_3 U_3>_0
                          -<\gamma^*_1 \gamma^*_2 \gamma_3>_0 <\gamma^*_3 U_2 U_3>_0 \big)
                          \big\} \nonumber \\
&&+<\gamma_1 \gamma^*_1 \gamma_2>_0 \big\{
          \partial_{U_1} \partial_{V_2} \big( <\gamma^*_2 V_1 \gamma_3>_0 <\gamma^*_3 U_2 U_3>_0
                                             -<\gamma^*_2 U_3 \gamma_3>_0 <\gamma^*_3 U_2 V_1>_0 \big) \nonumber \\
&&\hspace{90pt} +  \partial_{U_3} \partial_{V_2} \big( <\gamma^*_2 U_1 \gamma_3>_0 <\gamma^*_3 U_2 V_1>_0
                                          -<\gamma^*_2 U_2 \gamma_3>_0 <\gamma^*_3 U_1 V_1>_0 \big) \big\}. \nonumber
\end{eqnarray}
The right hand side of this equation is $0$ by equation \eqref{eqn:WDVVaP}.
The proposition is thus proved.
$\Box$

By Theorem \ref{thm:m=2} and Propositions \ref{B21=0} and \ref{B111=0}, we finish the proof of Theorem \ref{main theorem}.


\vspace{30pt} \noindent
Xiaobo Liu \\
School of Mathematical Sciences \& \\
Beijing International Center for Mathematical Research, \\
Peking University, Beijing, China. \\
Email: {\it xbliu@math.pku.edu.cn}
\ \\ \ \\
Chongyu Wang \\
School of Mathematical Sciences, \\
Peking University, Beijing, China. \\
Email: {\it wwccyy@pku.edu.cn}

\end{document}